\documentclass[twocolumn,letterpaper]{IEEEAerospaceCLS}  

\usepackage[]{graphicx}    
\newcommand{\ignore}[1]{}  
\usepackage{cite}
\usepackage{amsmath}
\usepackage{amssymb}
\usepackage{comment}
\usepackage{algorithm}
\usepackage{algpseudocode}
\usepackage{nicefrac, xfrac}

\begin{document}
\title{Fast Fuel-Optimal Constrained Impulsive Control with Application to Distributed Spacecraft}

\thispagestyle{plain}
\pagestyle{plain}

\author{%
Matthew Hunter\\ 
Department of Aeronautics \& Astronautics\\
Stanford University\\
Stanford, C2A 94305\\
mhunter8@stanford.edu
\and 
Simone D'Amico\\
Department of Aeronautics \& Astronautics\\
Stanford University\\
Stanford, CA 94305\\
damicos@stanford.edu
\thanks{\footnotesize 979-8-3503-5597-0/25/$\$31.00$ \copyright2025 IEEE}              
}

\maketitle

\thispagestyle{plain}
\pagestyle{plain}

\begin{abstract}
This work extends a reachable set theory approach for fuel-optimal control to incorporate convex and non-convex state, time, and magnitude constraints. Reachable set theory considers problems with norm-like cost functions and time-varying linear dynamics over a fixed time horizon and can produce optimality conditions for control input profiles. These optimality conditions enable the identification of optimal maneuver times and directions without needing to explicitly solve for the maneuvers themselves. While previous works have leveraged this approach to generate unconstrained, numerically-efficient optimal control solutions, this work proposes a novel impulsive control solver that finds provably-optimal trajectories under non-convex constraints. The solver is designed from a new set of optimality conditions, derived by reformulating a single reachable set into multiple time-ordered, step-wise reachable set waypoint sub-problems, where a waypoint is the state reached after a given maneuver. Non-convex constraints are satisfied at each waypoint, preserving convexity at each sub-problem without needing to convexify the full control problem, enabling the computationally efficient identification of constrained optimal maneuvers, and guaranteeing a user-defined bound on constrained cost. To produce an optimal maneuver plan, the solver rolls out constrained trajectories until the remaining unconstrained solution adheres to the problem constraints. The fuel-optimal control problem considered by this work is broadly applicable to many fields and particularly significant to Distributed Space Systems (DSS), which incorporate multiple spacecraft to overcome the limitations and risks posed by a single spacecraft. Multi-agent guidance and propulsive control techniques, often for small spacecraft at close separations, require guarantees of collision avoidance, efficient maneuvering to preserve a restricted fuel payload, and numerically-simple algorithms functional under on-board computational processing limitations. The proposed solver offers both user-defined bounds on constrained maneuver optimality and enhanced computational efficiency to directly address these challenges, and these capabilities are demonstrated comparatively against baseline solvers from the literature through Monte Carlo analysis on a variety of relevant DSS scenarios. Overall, this work enables the benefits of reachable set theory, namely computational efficiency and provable optimality, to be extended to control problems of higher fidelity.
\end{abstract}

\tableofcontents

\section{Introduction}

The fuel-optimal maneuver planning of constrained, discrete systems is significant and broadly applicable to many fields where the dynamics can be modelled as a Linear Time-Variant (LTV) system. One such important use case is propulsive spacecraft control, where cold-gas or combustive thrusts are approximately impulsive, or short in duration with respect to orbital period, and orbital motion is accurately characterized by linear models. In this field and this work, the term "fuel-optimal" refers to the minimum quantity of an actuation resource on-board a spacecraft, typically propellant for a thruster, required to accomplish a desired change in state over a specified time period. The fuel-optimal control problem is of particular interest to Distributed Space Systems (DSS), which incorporate multiple spacecraft to overcome the limitations and risks posed by a single spacecraft. Multi-agent guidance and control techniques, often for small spacecraft at close separations, require guarantees of collision avoidance, efficient maneuvering to preserve a restricted fuel payload, and numerically-simple algorithms functional under on-board computational processing limitations. Furthermore, the ability to solve more-complex LTV control problems can enhance the iterative control of higher-fidelity nonlinear dynamic systems, through techniques such as Sequential Convex Programming or an Extended Kalman Filter.

This class of control problem has received extensive focus in previous works, using a wide variety of solution approaches. The bulk of the literature can be divided into three categories: closed-form solutions, direct optimization, and indirect optimization. Closed-form solutions are commonly used for autonomous on-board maneuver planning, as they are often robust against dynamic perturbations, easy to test for initial condition sensitivity, and simple enough to be run on-board computationally limited hardware. An extensive DSS flight heritage has demonstrated these qualities on a variety of missions, including PRISMA and TanDEM-X \cite{damico_phd}. However, these solutions, and their optimality, are primarily specific to the state representation, dynamics model, and cost function of the use case for which the solutions are derived and often can only include limited perturbations. Direct optimization techniques are more generalizable to a variety of control problems, searching over maneuver times, directions, and magnitudes to evaluate the resulting effect on the overall control problem cost \cite{kim,sobiesiak,gong}. However, this search is computationally expensive, performed without insight into maneuver optimality, and does not commonly guarantee global optimality, as the control problem is usually non-convex and solved over a limited time horizon. In contrast, indirect optimization strategies define the necessary and sufficient conditions of optimality for the maneuver variables, primarily through convex optimization gradient techniques, to find a global provably-optimal solution to the control problem. The numerical efficiency of any given indirect optimizer depends heavily on the number of control variables and optimality conditions identified for the control problem.

For these reasons, indirect methods have been the most popular strategy for trajectory optimization in the literature. Approaches vary within this category depending on the type of problem considered: unconstrained, only considering the system dynamics, or constrained, incorporating operational constraints into the control problem. For unconstrained problems, previous works over the past several decades have leveraged both primer-dual theory and reachable set theory as optimization methodologies. Primer-dual theory uses a form of Lawden's "primer vector" to dictate control input through Pontryagin’s maximum principle and poses optimal control as a two-point boundary value problem that must satisfy a set of analytical conditions on the primer vector \cite{lawden}. Roscoe et al. created a trajectory optimizer for eccentric perturbed orbits using this strategy, but it required an accurate initialization of optimal maneuver timing and solved iteratively with significant computational load \cite{roscoe}. Reachable set theory analyzes the set of all possible maneuvers within the specified control period to find the most efficient maneuvers in the direction of the desired final state. Gilbert et al. designed a reachable set theory algorithm for problems with norm-like cost functions to converge to a globally optimal maneuver plan \cite{gilbert}. Koenig et al. extended this approach to time-variant cost functions, including periods of the control window where control actions are not allowed \cite{koenig_RST}. It's worth noting that the geometric optimality insights offered by reachable set theory have also enabled the creation of provably-optimal closed-form solutions for a variety of unconstrained DSS use-cases \cite{chernick,hunter_hybrid,hunter_robust}.

For operationally-constrained problems, previous authors have drawn from both types of numerical approaches to satisfy convex and non-convex constraints but have not been able to guarantee global optimality without significant simplifications. Sarno et al. balanced safe DSS trajectories with propellant consumption by using a centralized genetic algorithm to solve path planning and decentralized convex optimization to find individual cluster satellite minimum cost trajectories \cite{sarno}. Sequential Convex Programming (SCP) iterates between a higher fidelity non-convex control problem and lower fidelity convexified control problem to converge to a trajectory that satisfies non-convex mission constraints at efficient cost. Boelli et al. and Armellin have implemented SCP to guarantee passive and active collision avoidance, the former using Relative Orbital Elements (ROE) in near-circular orbit to inform definitions of safety and the latter using a cartesian state representation for orbits of arbitrary eccentricity \cite{borelli,armellin}. Guffanti et al. expanded these capabilities by accommodating thrust constraints and using Variation of Parameters to be robust against navigation uncertainty \cite{guffanti}. Model Predictive Control (MPC) iteratively finds a locally optimal control profile over a receding time horizon that satisfies constraints through a weighted objective function, and several authors have used convex, differentiable objective functions to reject navigation uncertainty, maintain safe inter-spacecraft separations, and remain within allowable thrust magnitude bounds \cite{patel,smith,sasaki}. Pseudo-spectral methods have been applied to large-scale DSS swarms, convexifying non-convex minimum separation constraints into nearest tangent plane constraints and breaking computationally-expensive full swarm collision-free guidance into decentralized sub-problems \cite{chen1,chen2}. Reachable set theory has been used in the space of constrained control as well. Cui et al. transformed reachable set-based ellipsiodal collision constraints into affine constraints for a computationally efficient SCP, while Bowerfind et al. rapidly estimated reachable sets by sampling minimum-fuel trajectories to optimize low-thrust rendezvous time \cite{cui,bowerfind}. Other optimization approaches for constrained control include Riemannian manifold path planning, particle swarm tree search, and reinforcement learning, but these provide no insight into the nature of optimal maneuvering \cite{kaptui,deka,xu}. Furthermore, none of these approaches are able to guarantee or provide a metric of global optimality for the full control problem without simplified, approximate constraints or substantial computational load.

In contrast, this work proposes a novel impulsive control solver that finds provably-optimal trajectories under both convex and non-convex constraints by extending a reachable set theory approach for fuel-optimal control. In this context, reachable set theory considers problems with norm-like cost functions and time-varying linear dynamics over a fixed time horizon and can produce optimality conditions for control input profiles. These optimality conditions enable the computationally efficient identification of optimal maneuver times and directions without needing to explicitly solve for the maneuvers themselves. The proposed solver is designed from a new set of optimality conditions, derived from reachable sets equivalently formulated as both time-ordered and as the convex combination of multiple reachable sets of equal combined cost. This results in a step-wise sum of single-maneuver waypoint sub-problems, where a waypoint is the state reached and time after a given maneuver. Convex and nonconvex state, time, and magnitude constraints are satisfied at each waypoint, preserving convexity without needing to convexify the full control problem. The proposed solver uses these optimality conditions to explore optimal maneuver branches from each waypoint, rolling out maneuver magnitudes and corresponding trajectories forward in time from the initial state and backward in time from the desired state to produce a maneuver plan that satisfies constraints with bounded optimality and limited computational load. 

Overall, this paper makes four major contributions to the state of the art, as detailed in four dedicated sections. First, the definition of a reachable set used by prior works is expressed as the intersection of multiple reachable sets, such that the original optimal control problem can be broken into optimal sequential sub-problems of the same form. This allows new necessary and sufficient optimality conditions to be defined relative to each sub-problem and enables the benefits of reachable set theory to be extended to control solvers that rely on optimal sub-structure. Second, constraints are defined and satisfied at the sub-problem level of the reachable set theory sequential formulation. Novel necessary and sufficient conditions for constrained optimality are derived with respect to maneuver magnitude within each sub-problem, such that trajectory optimization is solved through an optimal maneuver magnitude search. Third, the optimality conditions are leveraged to create an impulsive maneuver planner that balances computational efficiency and optimality to produce fast solutions to the constrained control problem with guaranteed optimality bounds. Fourth, potential DSS applications are demonstrated in simulation. Solver optimality, computational efficiency, and robustness to initial conditions are validated comparatively with baseline solvers from the literature on example DSS reconfigurations.

\section{Background}
\subsection{Problem Definition}

This paper addresses the impulsive or instantaneous control of linear time-varying dynamic systems, as defined by the following. Consider an LTV dynamics system, given as
\begin{equation}
    \boldsymbol{x}_f = \boldsymbol{\Phi}(t_0)\boldsymbol{x}_0 + \int_{t_0}^{t_f} \boldsymbol{\Phi}(\tau)\boldsymbol{B}(\tau)\boldsymbol{u}(\tau)d\tau,
\end{equation}
where $\boldsymbol{x}_0 \in \Re^n$ is the initial state at time $t_0$, $\boldsymbol{x}_f \in \Re^n$ is the final state at time $t_f$, $\boldsymbol{\Phi}(t) \in \Re^{n \times n}$ is the State Transition matrix that propagates a given state from time $t$ to $t_f$, and $\boldsymbol{B}(t) \in \Re^{n \times m}$ is the control matrix that converts control actions $\boldsymbol{u}(t) \in \Re^m$ to changes in state $\boldsymbol{x}$. For simplicity, the pseudostate $\boldsymbol{\omega} \in \Re^n$, the dynamic invariants of the control problem, and control matrix $\boldsymbol{\Gamma}(t) \in \Re^{n \times m}$, the translation from control actions to changes in final state, are defined as,
\begin{equation}
    \boldsymbol{\omega} = \boldsymbol{x}_f - \boldsymbol{\Phi}(t_0)\boldsymbol{x}_0, \hspace{5mm} \boldsymbol{\Gamma}(t) = \boldsymbol{\Phi}(t)\boldsymbol{B}(t).
\end{equation}
Therefore, the unconstrained optimal control problem considered by this paper is defined as
\begin{align}
    \text{minimize: }& \int_{t_0}^{t_f} f(\boldsymbol{u}(\tau),\tau)d\tau \label{eq: continuous problem} \\
    \text{subject to: }& \boldsymbol{\omega} = \int_{t_0}^{t_f} \boldsymbol{\Gamma}(\tau)\boldsymbol{u}(\tau)d\tau. \nonumber
\end{align}
Again, the term "unconstrained" is used here to refer to control problems that are only constrained by the system dynamics. Cost function $f$ is piecewise-defined as a set of norm-like functions $f_j$ over an arbitrary number $o$ of mutually exclusive time intervals, defined as
\begin{equation}
    f(\boldsymbol{u}(t),t) = 
    \begin{cases}
        f_1(\boldsymbol{u}(t)), & t_{0,1} \leq t \leq t_{f,1} \\
        f_2(\boldsymbol{u}(t)), & t_{0,2} \leq t \leq t_{f,2} \\
         & ... \\
        f_o(\boldsymbol{u}(t)), & t_{0,o} \leq t \leq t_{f,o} \\
        0 & \text{otherwise}
    \end{cases}, \label{eq: cost function}
\end{equation}
where $t_{0,j}$ and $t_{f,j}$ are the beginning and end of each time interval. The norm-like functions $f_j$ satisfy three properties: 1) $f_j$ is defined for all maneuvers occurring within $t_{0,j} \leq t \leq t_{f,j}$, 2) all sublevel sets of $f_j(\boldsymbol{u}(t),t)$ are convex and compact, and 3) $f_j(\alpha\boldsymbol{u}(t),t) = \alpha f_j(\boldsymbol{u}(t),t)$ for any $\alpha \geq 0$. Examples of norm-like cost functions are given in Table \ref{tab: cost functions}, as well as by previous authors in the literature \cite{koenig_RST}.

\begin{table}[b]
    \centering
    \caption{\bf Example Cost Functions and Corresponding Spacecraft Propulsion Setup Constraints for $\boldsymbol{u}(t) \in \Re^3$}
    \vspace{-1mm}
    \begin{tabular}{c|c}
        Cost Function & Propulsion Constraint \\
        \hline \hline
        \rule{0pt}{2ex} $||\boldsymbol{u}||_2$ & Align a single thruster with \\
         & the desired maneuver direction. \\
        \hline
        \rule{0pt}{2ex}  & Fixed attitude with three pairs of \\
        $||\boldsymbol{u}||_1$ & thrusters mounted on opposite sides \\ 
         & on mutually perpendicular axes. \\
        \hline
        \rule{0pt}{2ex} & Two pairs of thrusters on \\
        $|u_1| + \sqrt{u_2^2 + u_3^2}$ & perpendicular axes where one \\
         & axis is fixed. \\
         \hline
    \end{tabular}
    \label{tab: cost functions}
\end{table}

This work focuses on impulsive control actions, or instantaneous changes in state, assuming that the time required to complete each maneuver is short with respect to the secular and long-term periodic effects present in $\boldsymbol{\Phi}(t)$. Given there are no restrictive magnitude constraints on $\boldsymbol{u}$, time-varying control actions can be expressed as a set of $k$ impulses, given as
\begin{equation}
	\boldsymbol{u} = \sum_{j=1}^k \boldsymbol{u}_j \delta (t-t_j),
\end{equation}
where $\delta$ denotes the Dirac delta function. Using the impulsive assumption, the unconstrained control problem considered by the remainder of this paper is given as
\begin{equation}
	\text{minimize: } c = \sum_{j=1}^{k} f_j(\boldsymbol{u}_j,t_j) \hspace{3mm} \text{subject to: } \boldsymbol{\omega} = \sum_{j=1}^{k} \boldsymbol{\Gamma}(t_j)\boldsymbol{u}_j, \label{eq: impulsive problem}
\end{equation}
where all $t_j \in T$ and $c$ is the minimum cost objective.

\subsection{Reachable Set Theory}

Previous applications of reachable set theory for optimal control have focused on time-varying dynamical systems without operational constraints, where all control actions are permissible independent of state, time, and maneuver order. The analysis summarized here forms the foundation of the advances made in this work and will later be built upon to incorporate constraints beyond the LTV dynamics. Let $U(c,t)$ be the set of control actions $\boldsymbol{u}$ with a cost no greater than $c$ at time $t$, defined as
\begin{equation}
    U(c,t) = \{\boldsymbol{u}: f(\boldsymbol{u},t) \leq c \}. \label{eq: U}
\end{equation}
Let $S(c,t)$ be the set of all $\boldsymbol{\omega}$ that can be reached by a single maneuver $\boldsymbol{u}$ of cost no greater than $c$ at time $t$, given as
\begin{equation}
    S(c,t) = \{\boldsymbol{y}: \boldsymbol{y} = \boldsymbol{\Gamma}(t)\boldsymbol{u}, \boldsymbol{u} \in U(c,t) \}. \label{eq: S(c,t)}
\end{equation}
This set is not generally convex. Let $T$ be the discontinuous set of times $t$ for admissible control actions in a control window bounded by $t_0$ and $t_f$, given as
\begin{equation}
    T = T(t_0,t_f) = \{t: t_0 \leq t \leq t_f \}. \label{eq: T}
\end{equation}
Then, the set of all $\boldsymbol{\omega}$ that can be reached by a single maneuver $\boldsymbol{u}$ of cost no greater than $c$ at any any time in $T$ is defined as
\begin{equation}
    S(c,T) = \cup_{t \in T} S(c,t), \label{eq: S(c,T)}
\end{equation}
which is also generally non-convex. Finally, let $S^*(c,T)$ be the set of $\boldsymbol{\omega}$ that can be reached by $k \geq 1$ control actions of combined cost no greater than $c$ at times in $T$, defined equivalently as
\begin{align}
    S^*(c,T) = \{\boldsymbol{z}: \boldsymbol{z} &= \sum\nolimits_{j=1}^{k} \boldsymbol{\Gamma}(t_j)\boldsymbol{u}_j, t_j \in T, \label{eq: S* base} \\
    &\boldsymbol{u}_j \in U(c_j,t_j), \sum\nolimits_{j=1}^{k} c_j = c \} \nonumber \\
    = \{\boldsymbol{z}: \boldsymbol{z} &= \sum\nolimits_{j=1}^{k}\alpha_j\boldsymbol{y}_j, \boldsymbol{y}_j \in S(c,T), \label{eq: S* convex hull} \\
    &\alpha_j \geq 0, \sum\nolimits_{j=1}^{k} \alpha_j = 1 \}. \nonumber
\end{align}
This shows that $S^*(c,T)$ is a linear combination and convex hull of $\boldsymbol{\omega}$ in $S(c,T)$. 

Returning to the optimal control problem in (\ref{eq: impulsive problem}), the set $S^*(c,T)$ in (\ref{eq: S* convex hull}) completely describes the dynamics constraint such that the problem can be defined as
\begin{equation}
    \text{minimize: } c \hspace{5mm} \text{subject to: } \boldsymbol{\omega} \in S^*(c,T). \label{eq: problem with RST}
\end{equation}
Previous authors have leveraged two important qualities of this formulation to produce computationally efficient and intuitive optimal control solvers. First, the norm-like cost function $f$ and pseudostates $\boldsymbol{\omega} \in S^*(c,T)$ scale linearly with magnitude. Second, the fact that $S^*(c,T)$ is a convex set means that it must be able to be described by a closed intersection of hyperplanes. Therefore, the scalar $c$ that scales $S^*(c,T)$ such that its boundary contains $\boldsymbol{\omega}$ is the minimum or optimal cost of (\ref{eq: problem with RST}). This interpretation is realized as the dual of (\ref{eq: problem with RST}), defined as 
\begin{equation}
    \text{maximize: } c \hspace{5mm} \text{subject to: } \max_{\boldsymbol{z} \in S^*(c,T)} \boldsymbol{\lambda}^T\boldsymbol{z} \leq \boldsymbol{\lambda}^T\boldsymbol{\omega}, \label{eq: RST dual}
\end{equation}
where $\boldsymbol{\lambda}$ is a supporting hyperplane of $S^*(c,T)$. This dual formulation gains further utility by expanding one side of the constraint using the linear relation between $S(c,T)$ and $S^*(c,T)$, given as
\begin{align}
    \max_{\boldsymbol{z} \in S^*(c,T)} \boldsymbol{\lambda}^T\boldsymbol{z} &= \max_{\boldsymbol{y} \in S(c,T)} \boldsymbol{\lambda}^T\boldsymbol{y} \label{eq: contact function} \\
    & = c \left( \max_{t \in T} \left( \max_{\boldsymbol{u} \in U(1,t)} \boldsymbol{\lambda}^T\boldsymbol{\Gamma}(t)\boldsymbol{u} \right) \right). \nonumber 
\end{align}
Defining a contact function generally as
\begin{equation}
    g_Z(\boldsymbol{\lambda}) = \max_{\boldsymbol{z} \in Z} \boldsymbol{\lambda}^T\boldsymbol{z},
\end{equation}
for an arbitrary set $Z$, the constraint expansion allows (\ref{eq: RST dual}) to be reformulated as Second-Order Cone Program (SOCP), given as
\begin{equation}
    \text{maximize: } c \hspace{5mm} \text{subject to: } c \leq \frac{\boldsymbol{\lambda}^T\boldsymbol{\omega}}{\max_{t \in T} g_{S(1,t)}(\boldsymbol{\lambda})}, \label{eq: RST SOCP cost}
\end{equation}
or more simply after algebraic reformulation as
\begin{equation}
    \text{maximize: } \boldsymbol{\lambda}^T\boldsymbol{\omega} \hspace{5mm} \text{subject to: } \max_{t \in T} g_{U(1,t)}(\boldsymbol{\Gamma}^T(t) \boldsymbol{\lambda}) \leq 1. \label{eq: RST SOCP}
\end{equation}

\subsection{Optimality Conditions}

\begin{figure} [b]
    \centering
    \includegraphics[width=1\linewidth]{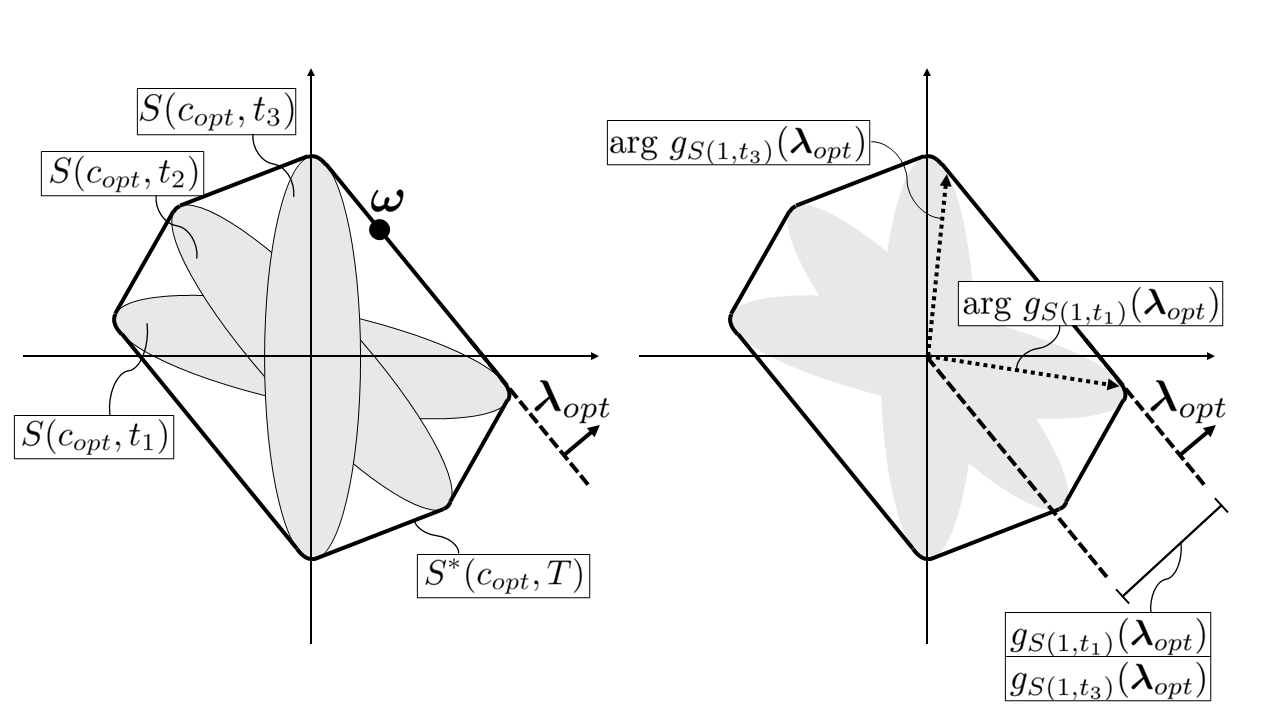}
    \caption{2-dimensional reachable set formed as convex hull of states reached by single maneuvers at times $t_1 < t_2 < t_3$ (left). Optimality conditions evaluated on the boundary of $S^*(1,T)$ through contact functions (right). }
    \label{fig: unconstrained}
\end{figure}

Previous authors have found optimality conditions for unconstrained control through (\ref{eq: RST SOCP}), given that all control actions can be analyzed a-priori before solving any part of the maneuver plan. These conditions can be summarized as follows. First, consider an optimal solution to (\ref{eq: RST SOCP}) defined as $c_{opt}$ and $\boldsymbol{\lambda}_{opt}$. $c_{opt}$ is the minimum cost to reach $\boldsymbol{\omega}$ if and only if $\boldsymbol{\omega}$ is on the boundary of $S^*(c_{opt},T)$. This boundary is described by the supporting hyperplane of $S^*(c_{opt},T)$ defined perpendicular to $\boldsymbol{\lambda}_{opt}$, and the optimal control input profile may not be unique if $\boldsymbol{\omega}$ lies on an edge or point of the boundary of $S^*(c_{opt},T)$. Second, for any $\boldsymbol{\omega}$ on the boundary of $S^*(c_{opt},T)$, an optimal control input sequence must exist such that $k \leq n$. Therefore, an optimal control sequence must exist for any $\boldsymbol{\omega}$ of no more maneuvers than the dimensionality of $\boldsymbol{\omega}$. Third, control actions should only occur when $g_{U(1,t)}(\boldsymbol{\Gamma}^T(t) \boldsymbol{\lambda}_{opt})$ reaches its maximum value in the control window $T$. This condition is reflected in the constraint of (\ref{eq: RST SOCP}). Fourth, control actions should be in the direction of $\text{arg }g_{U(1,t)}(\boldsymbol{\Gamma}^T(t)\boldsymbol{\lambda}_{opt})$. The third and fourth conditions can be evaluated for an arbitrary $\boldsymbol{\lambda}$ and time $t$ for an L2-norm cost function through the Cauchy Schwartz inequality, given as
\begin{align}
    g_{U(1,t)}(\boldsymbol{\Gamma}^T(t) \boldsymbol{\lambda}) &= ||\boldsymbol{\Gamma}^T(t) \boldsymbol{\lambda}||_2, \label{eq: g_U} \\
    \text{arg }g_{U(1,t)}(\boldsymbol{\Gamma}^T(t) \boldsymbol{\lambda}) &= \frac{\boldsymbol{\Gamma}^T(t) \boldsymbol{\lambda}}{||\boldsymbol{\Gamma}^T(t) \boldsymbol{\lambda}||_2}. \label{eq: arg g_U}
\end{align}
These optimality conditions are visualized in Figure \ref{fig: unconstrained}.

\section{Sequential Formulation}
\subsection{Intersecting Reachable Sets}

Introducing time and state constraints to the control problem violates the assumption that all maneuvers are permissible independent of time, state, and maneuver order. Therefore, the definitions of the reachable sets shall be reformulated such that they are defined, and can be solved, as sequential reachable set sub-problems. First, the maneuvers within $S^*(c,T)$ in (\ref{eq: S* base}) are split into $l$ groups, given as
\begin{align}
    S^*(c,T) = \{\boldsymbol{z}: \boldsymbol{z} &= \sum\nolimits_{i=1}^l \sum\nolimits_{j=1}^{k_i} \boldsymbol{\Gamma}(t_{i,j})\boldsymbol{u}_{i,j}, \label{eq: S* grouped} \\
    & t_{i,j} \in T, \boldsymbol{u}_{i,j} \in U(c_{i,j},t_{i,j}), \nonumber \\
    & \sum\nolimits_{i=1}^l \sum\nolimits_{j=1}^{k_i} c_{i,j} = c, \sum\nolimits_{i=1}^l k_i = k \}, \nonumber 
\end{align}
where $k_i$ is the number of maneuvers within group $i$. Consider time-ordering each group $i$, given as
\begin{align}
    S^*(c,T) = \{\boldsymbol{z}: \boldsymbol{z} &= \sum\nolimits_{i=1}^l \sum\nolimits_{j=1}^{k_i} \boldsymbol{\Gamma}(t_{i,j})\boldsymbol{u}_{i,j},  \\
    & t_{i,j} \in T, t_{i,j} \geq t_{i-1,j}, \boldsymbol{u}_{i,j} \in U(c_{i,j},t_{i,j}), \nonumber \\
    &\sum\nolimits_{i=1}^l \sum\nolimits_{j=1}^{k_i} c_{i,j} = c, \sum\nolimits_{i=1}^l k_i = k \}. \nonumber
\end{align}
This constraint may be added without loss of generality, as the ordering has no effect on the shape of $S^*(c,T)$ or the resulting $\boldsymbol{\lambda}_{opt}$ from (\ref{eq: RST SOCP}) that defines the optimality of a given maneuver time and direction. This enables $S^*(c,T)$ to be rearranged to an equivalent time-ordered step-wise formulation, given as
\begin{align}
    S^*(c&,T) = \{\boldsymbol{z}: \boldsymbol{z} = \boldsymbol{z}_l, \boldsymbol{z}_i = \sum\nolimits_{j=1}^{k_i} \boldsymbol{\Gamma}(t_{i,j}) \boldsymbol{u}_{i,j} + \boldsymbol{z}_{i-1}, \nonumber \\
    &\boldsymbol{z}_{0} = \Vec{\boldsymbol{0}}, t_{i,j} \in T, t_{i,j} \geq t_{i-1,j}, \boldsymbol{u}_{i,j} \in U(c_{i,j},t_{i,j}), \nonumber \\
    &\sum\nolimits_{i=1}^l \sum\nolimits_{j=1}^{k_i} c_{i,j} = c, \sum\nolimits_{i=1}^l k_i = k \}, \label{eq: S* step-wise}
\end{align}
where each $\boldsymbol{z}_i$ is a waypoint of the maneuver plan, defined as an intermediate pseudostate achieved by a set of maneuvers in $U(c_{i,j},t_{i,j})$. The time-ordering in (\ref{eq: S* step-wise}) allows the waypoints to describe the trajectory of a maneuver plan across pseudostate-space. It is clear from this formulation that the maneuvers between each waypoint match the formulation of $S^*(c,T)$ in (\ref{eq: S* base}), describing a convex hull of size $\sum_{j=1}^{k_i} c_{i,j}$. Therefore, consider subsets of the control window $T$ defined from a particular time $t_j$ in the control window as
\begin{align}
    T^+_j &= T^+(t^+_j) = \{t: t \geq t^+_j, t \in T \}, \label{eq: T+} \\
    T^-_j &= T^-(t^-_j) = \{t: t \leq t^-_j, t \in T \}, \label{eq: T-}
\end{align}
where the subset after $t_j$ is denoted by superscript "+" and the subset before $t_j$ is denoted by superscript "-". Incorporating (\ref{eq: T+}) and (\ref{eq: T-}) into (\ref{eq: S* step-wise}), $S^*(c,T)$ can be reformulated as
\begin{align}
    S^*(c,T) = \{\boldsymbol{z}:& \boldsymbol{z} = \sum\nolimits^l_{i=1} \boldsymbol{z}^\Delta_{i}, \boldsymbol{z}^\Delta_{i} \in S^*(c_i,T_i^+), \label{eq: S* forward} \\
    &T^+_i \in T^+_{i-1}, \sum\nolimits_{i=1}^l c_i = c \}, \nonumber
\end{align}
or equivalently backward in time as
\begin{align}
    S^*(c,T) = \{\boldsymbol{z}:& \boldsymbol{z} = \sum\nolimits^l_{i=1} \boldsymbol{z}^\Delta_{i}, \boldsymbol{z}^\Delta_{i} \in S^*(c_i,T_i^-), \label{eq: S* backward} \\
    &T^-_i \in T^-_{i-1}, \sum\nolimits_{i=1}^l c_i = c \}, \nonumber
\end{align}
where $\boldsymbol{z}^\Delta_{i}$ is the change in pseudostate from waypoint $\boldsymbol{z}_{i-1}$ to waypoint $\boldsymbol{z}_i$. $S^*(c,T)$ is now posed as a sequential time-ordered series of waypoint reachable set sub-problems. It is worth noting here that all time-ordered formulations in this paper can be equivalently posed either forward in time, from $t_0$ to $t_f$, or backward in time, from $t_f$ to $t_0$, with very similar formulations. Therefore for simplicity, both directions will be expressed concisely as
\begin{align}
    S^*(c,T) = \{\boldsymbol{z}:& \boldsymbol{z} = \sum\nolimits^l_{i=1} \boldsymbol{z}^\Delta_{i}, \boldsymbol{z}^\Delta_{i} \in S^*(c_i,T_i^{(\cdot)}), \label{eq: S* time-ordered} \\
    &T^{(\cdot)}_i \in T^{(\cdot)}_{i-1}, \sum\nolimits_{i=1}^l c_i = c \}, \nonumber
\end{align}
where $T^{(\cdot)}_i = T^+_i$ or $T^{(\cdot)}_i = T^-_i$ depending on the chosen direction of time-ordering.

Using the new sequential set definitions, the optimal control problem can now also be posed sequentially. Consider a sequential expansion of the optimal problem in (\ref{eq: RST dual}), given as
\begin{align}
    \text{maximize: }& \sum_{i=1}^l c_i \label{eq: RST sequential 1} \\
    \text{subject to: }& \sum_{i=1}^l \max_{\boldsymbol{z}^\Delta_i \in S^*(c_i,T_i^{(\cdot)})} \boldsymbol{\lambda}_i^T\boldsymbol{z}^\Delta_i \leq \sum_{i=1}^l \boldsymbol{\lambda}_i^T\boldsymbol{\omega}_i,  \nonumber
\end{align}
where $\sum_{i=1}^l \boldsymbol{\lambda}^T_i\boldsymbol{\omega}_i = \boldsymbol{\lambda}^T\boldsymbol{\omega}$. Recalling that this sum of convex hulls is itself a convex hull, this sequential expansion can be expressed as
\begin{align}
    \text{maximize: }& \sum_{i=1}^l \beta_i c \label{eq: RST sequential 2} \\
    \text{subject to: }& \sum_{i=1}^l \max_{\boldsymbol{z}^\Delta_i \in S^*(\beta_i c,T_i^{(\cdot)})} \boldsymbol{\lambda}^T\boldsymbol{z}^\Delta_i \leq \sum_{i=1}^l \beta_i\boldsymbol{\lambda}^T\boldsymbol{\omega}, \nonumber
\end{align}
where $\beta_i \geq 0$ and $\sum_{i=1}^l \beta_i = 1$. Consider breaking this problem into a series of sequential sub-problems, each sub-problem $i$ given as
\begin{align}
    \text{maximize: }& \beta_i c \label{eq: RST sequential sub-problem} \\
    \text{subject to: }& \max_{\boldsymbol{z}^\Delta_i \in S^*(\beta_i c,T_i^{(\cdot)})} \boldsymbol{\lambda}^T\boldsymbol{z}^\Delta_i \leq  \beta_i\boldsymbol{\lambda}^T\boldsymbol{\omega}. \nonumber
\end{align}
It can be seen by visual inspection that this sub-problem is almost identical in structure to the full problem in (\ref{eq: RST dual}), such that $\beta_i = 1$ and $T_i^{(\cdot)} = T$ for any $i$ recovers the full problem definition in (\ref{eq: RST dual}). This reveals an optimal, overlapping sub-problem structure compatible with the principle of optimality commonly leveraged by dynamic programming, meaning that the optimal solution to (\ref{eq: RST dual}) must be comprised of the optimal solutions to its sub-problems \cite{kochenderfer}. Therefore, combining the optimal solutions to each sub-problem $i$ in (\ref{eq: RST sequential sub-problem}) results in an optimal solution to the full problem. Applying the algebraic manipulations used to reach (\ref{eq: RST SOCP}), each of these sub-problems can be formulated as an SOCP, given as
\begin{align}
    \text{maximize: }& \beta_i \boldsymbol{\lambda}^T\boldsymbol{\omega} \label{eq: RST SOCP sub-problem} \\
    \text{subject to: }&  \max_{t \in T_i^{(\cdot)}} g_{U(1,t)}(\boldsymbol{\Gamma}^T(t) \boldsymbol{\lambda}) \leq 1. \nonumber
\end{align}
It is important to note that $\beta_i$ only scales the objective function of this SOCP
and has no effect on $\boldsymbol{\lambda}_{opt}$, such that the solution $\boldsymbol{\lambda}_{opt}$ for each sub-problem matches the solution $\boldsymbol{\lambda}_{opt}$ to the full problem SOCP in (\ref{eq: RST SOCP}). However, this means that (\ref{eq: RST SOCP sub-problem}) cannot generate a meaningful solution for $\beta_i$, and this optimal value must be solved separately. 

\subsection{Optimality Conditions}

\begin{figure}[t]
    \centering
    \includegraphics[width=1\linewidth]{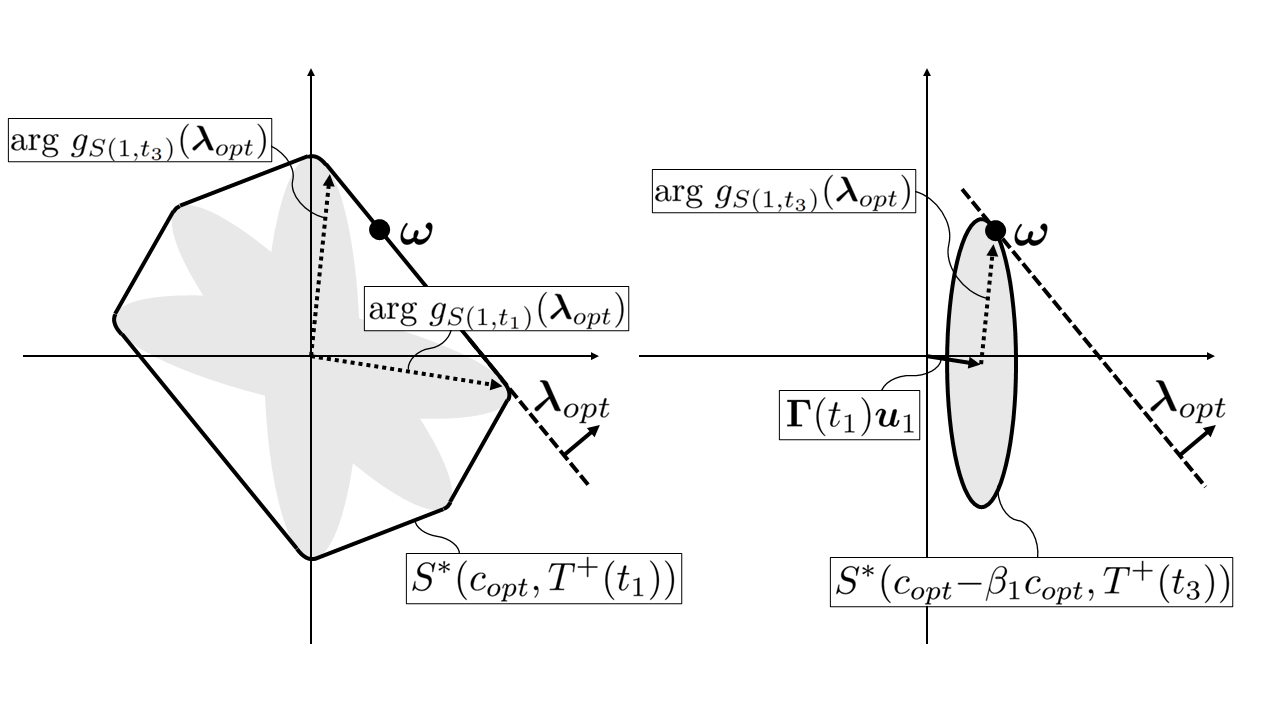}
    \caption{Reachable set forward time-ordered optimal sub-problems at $T^+(t_1)$ (left) and $T^+(t_3)$ (right) generated from the full problem reachable set in Figure \ref{fig: unconstrained}.}
    \label{fig: sequential}
\end{figure}

Given the close similarities between the full problem SOCP in (\ref{eq: RST SOCP}) and the sub-problem SOCP in (\ref{eq: RST SOCP sub-problem}), the necessary and sufficient conditions previously derived for optimal maneuver direction and timing are also applicable to (\ref{eq: RST SOCP sub-problem}). However, additional conditions are required to characterize an optimal $\beta_i$ for a given sub-problem $i$, specifically defining the transition between sub-problem $i$ and sub-problem $i+1$. To identify these conditions, consider splitting $\boldsymbol{\omega}$ into its sub-problem components, given as
\begin{align}
    \text{maximize: }& \beta_i \boldsymbol{\lambda}^T \left( \sum\nolimits_{i'=1}^{l} \boldsymbol{\omega}_{i'} \right) \label{eq: RST SOCP sub-problem split}\\
    \text{subject to: }&  \max_{t \in T_i^{(\cdot)}} g_{U(1,t)}(\boldsymbol{\Gamma}^T(t) \boldsymbol{\lambda}) \leq 1, \nonumber
\end{align}
where the primed indexing $i'$ is used to split $\boldsymbol{\omega}$ within the sub-problem already denoted by $i$. Assuming all previous time-ordered sub-problems were solved optimally, the optimal solution to (\ref{eq: RST SOCP sub-problem split}) must be equivalent to the reduced problem
\begin{align}
    \text{maximize: }& \beta_i \boldsymbol{\lambda}^T \Tilde{\boldsymbol{\omega}} \label{eq: RST SOCP sub-problem reduced} \\
    \text{subject to: }&  \max_{t \in T_i^{(\cdot)}} g_{U(1,t)}(\boldsymbol{\Gamma}^T(t) \boldsymbol{\lambda}) \leq 1, \nonumber
\end{align}
where $\Tilde{\boldsymbol{\omega}} = \sum\nolimits_{i'=i}^l \boldsymbol{\omega}_{i'}$. Next, consider restricting the set $S^*(\beta_i c, T_i^{(\cdot)})$ to a single maneuver. As discussed previously, the time-ordering constraint has no effect on optimality or the shape of $S^*(c,T)$, so the single maneuver assumption can be made without loss of generality. These properties are demonstrated by visual inspection in Figure \ref{fig: sequential} and detailed by the following. Consider a series of forward time-ordered sub-problems resulting from the full 2-dimensional problem in Figure \ref{fig: unconstrained}. $\boldsymbol{\lambda}_{opt}$ is known to be the solution to the full control problem in (\ref{eq: RST SOCP}) and as a result the solution to all sub-problems in (\ref{eq: RST SOCP sub-problem}). Optimal maneuvers are identified at $t_1$ and $t_3$ in the directions specified by (\ref{eq: arg g_U}). Consider restricting the analysis to a first sub-problem $i=1$ at $T^+(t_1)$ (Fig. \ref{fig: sequential}, left) and a second sub-problem $i=2$ at $T^+(t_3)$ (Fig. \ref{fig: sequential}, right). The first sub-problem can be equivalently constrained to maneuvers after $t_1$ rather than $t_0$ as no optimal maneuver occur between $t_0$ and $t_1$. The same argument can be made for the second sub-problem at $t_3$, resulting in the removal of $t_2$ from the second sub-problem. The first sub-problem cost $\beta_1c$ achieves its optimal value when the $\boldsymbol{\lambda}_{opt}$ solution to (\ref{eq: RST SOCP sub-problem reduced}) for the resulting second sub-problem at $T^+(t_3)$ matches the solution of (\ref{eq: RST SOCP sub-problem}), also resulting in the same optimal maneuver direction at $t_3$. 

Given an arbitrary unitary maneuver $\Hat{\boldsymbol{u}}_i$ at time $t_i \in T_i^{(\cdot)}$ at a given magnitude $\beta_i$, the resulting sub-problem $i+1$ is given as
\begin{align}
    \text{maximize: }& \beta_{i+1} \boldsymbol{\lambda}^T \left( \Tilde{\boldsymbol{\omega}} - \beta_i\boldsymbol{\Gamma}(t_i)\Hat{\boldsymbol{u}}_i \right) \label{eq: RST SOCP sub-problem i+1} \\
    \text{subject to: }&  \max_{t \in T_{i+1}^{(\cdot)}} g_{U(1,t)}(\boldsymbol{\Gamma}^T(t) \boldsymbol{\lambda}) \leq 1. \nonumber
\end{align}
This formulation enables necessary and sufficient optimality conditions to be established for $\beta_i$ with respect to the $\boldsymbol{\lambda}_{opt,i+1}$ that results as the solution to (\ref{eq: RST SOCP sub-problem i+1}), without needing to solve for $\beta_{i+1}$. First, $\Hat{\boldsymbol{u}}_{i}$ at time $t_{i}$ is optimal in (\ref{eq: RST SOCP sub-problem i+1}) if its contact function satisfies the constraint of (\ref{eq: RST SOCP sub-problem i+1}), given as
\begin{equation}
    g_{\Hat{\boldsymbol{u}}_{i}}(\boldsymbol{\Gamma}^T(t_{i})\boldsymbol{\lambda}_{opt,i+1}) \geq 1, \label{eq: g check}
\end{equation}
such that the contact function for $\Hat{\boldsymbol{u}}_{i}$ matches or exceeds the maximal value found for maneuvers over the subset of times in $ T_{i+1}^{(\cdot)}$. Second, recall that the solution $\boldsymbol{\lambda}_{opt,i+1}$ to (\ref{eq: RST SOCP sub-problem i+1}) yields the optimal maneuver times $t_{opt,i+1}$ and directions $\Hat{\boldsymbol{u}}_{opt,i+1}$ for an unconstrained reachable set, given as
\begin{align}
    t_{opt,i+1} &= \text{arg } \max_{t \in T_{i+1}^{(\cdot)}} g_{U(1,t)}(\boldsymbol{\Gamma}^T(t) \boldsymbol{\lambda}_{opt,i+1}), \\
    \Hat{\boldsymbol{u}}_{opt,i+1} &= \text{arg } g_{U(1,t_{opt,i+1})}(\boldsymbol{\Gamma}^T(t_{opt,i+1}) \boldsymbol{\lambda}_{opt,i+1}),
\end{align}
and evaluated through the Cauchy-Swartz inequality in (\ref{eq: g_U}) and (\ref{eq: arg g_U}). Then, $\Hat{\boldsymbol{u}}_{i}$ at time $t_{i}$ is optimal in (\ref{eq: RST SOCP sub-problem i+1}) if it acts in the same direction and with equal or greater efficiency in pseudostate space as an optimal maneuver identified in (\ref{eq: RST SOCP sub-problem i+1}), defined as
\begin{align}
    &\frac{\left(\boldsymbol{\Gamma}(t_{i})\Hat{\boldsymbol{u}}_{i}\right)^T\left(\boldsymbol{\Gamma}(t_{i+1,opt})\Hat{\boldsymbol{u}}_{i+1,opt}\right)}{||\boldsymbol{\Gamma}(t_{i})\Hat{\boldsymbol{u}}_{i}||_2 ||\boldsymbol{\Gamma}(t_{i+1,opt})\Hat{\boldsymbol{u}}_{i+1,opt}||_2} = 1 \hspace{2mm} \forall \Hat{\boldsymbol{u}}_{opt,i+1}, \label{eq: efficiency check 1} \\
    &\frac{\left(\boldsymbol{\Gamma}(t_{i})\Hat{\boldsymbol{u}}_{i}\right)^T\left(\boldsymbol{\Gamma}(t_{i+1,opt})\Hat{\boldsymbol{u}}_{i+1,opt}\right)}{||\boldsymbol{\Gamma}(t_{i+1,opt})\Hat{\boldsymbol{u}}_{i+1,opt}||^2_2} \geq 1 \hspace{2mm} \forall \Hat{\boldsymbol{u}}_{opt,i+1}, \label{eq: efficiency check 2}
\end{align}
for direction and efficiency respectively. Third, consider leveraging the definition of cost in the constraint of (\ref{eq: RST SOCP cost}) to realize the remaining cost of sub-problem $i$, given as
\begin{align}
    \text{maximize: }& \beta_i c \\
    \text{subject to: }& \beta_i c \leq \frac{\boldsymbol{\lambda}^T \Tilde{\boldsymbol{\omega}}}{\max_{t \in T_i^{(\cdot)}} g_{U(1,t)}(\boldsymbol{\Gamma}^T(t) \boldsymbol{\lambda})}. \nonumber
\end{align}
$\Hat{\boldsymbol{u}}_{i}$ at time $t_{i}$ is optimal in (\ref{eq: RST SOCP sub-problem i+1}) if the minimum cost of sub-problem $i+1$ is no greater than the minimum cost found in the previous sub-problem $i$ minus the cost of the maneuver, equal to $\beta_i$ due to the unitary $\Hat{\boldsymbol{u}}_{i}$, defined as
\begin{align}
    &\frac{\boldsymbol{\lambda}_{opt,i}^T\Tilde{\boldsymbol{\omega}}}{\max_{t \in T_{i}^{(\cdot)}} g_{U(1,t)}(\boldsymbol{\Gamma}^T(t) \boldsymbol{\lambda}_{opt,i})} - \beta_i \label{eq: min cost check} \\
    &= \frac{\boldsymbol{\lambda}_{opt,i+1}^T \left( \Tilde{\boldsymbol{\omega}} - \beta_i \boldsymbol{\Gamma}(t_{i})\Hat{\boldsymbol{u}}_{i} \right)}{\max_{t \in T_{i+1}^{(\cdot)}} g_{U(1,t)}(\boldsymbol{\Gamma}^T(t) \boldsymbol{\lambda}_{opt,i+1})}. \nonumber
\end{align}
It is worth noting here that the three necessary and sufficient optimality conditions presented here for $\beta_i$ are independent definitions of optimality. However, all three will be relied upon in the creation of a control solver to make the identification of the optimal $\beta_i$ robust against numerical instability. Trajectory optimization can now be accomplished by iterating between optimal maneuver time and direction identification through (\ref{eq: RST SOCP sub-problem}) and finding the optimal magnitude $\beta_i$ for a single unitary optimal maneuver at waypoint $i$.

\section{Constrained Reachable Sets}

This work will consider three types of constraints: state, time, and magnitude. State constraints define an area of pseudostate-space where waypoints are inadmissible and the trajectory between them cannot enter. State constraints can be non-convex, typically as a keep out exclusion zone, or convex, typically as a maximum deviation from a desired nominal state or as a halfspace. Time constraints define periods of the control window where control inputs are inadmissible. While previous authors have addressed convex time constraints, where the inadmissible control input times are known a-priori to solving the control problem, this work will address non-convex time constraints, where inadmissible control input times are defined relative to other maneuvers such that no two maneuvers may occur within a specified interval. Maneuver magnitude constraints are primarily expressed as an upper bound on the cost of a maneuver. These convex bounds are usually dictated by the hardware limitations of the actuation system, such as the max outflow of a thruster. To incorporate these constraints into the reachable set theory approach, the definitions of $S(c,T)$ and $S^*(c,T)$ will be expanded. Consider an arbitrary state constraint set $Q$, containing all states where waypoints and related pseudostate trajectories are inadmissible, a minimum maneuver time spacing constraint $\Delta t$, and a maximum maneuver cost constraint $c_{max}$. Beginning with constraint $Q$, let $S(c,t,Q,\boldsymbol{y}_0)$ be the set of all $\boldsymbol{\omega}$ in $S(c,t)$ that are not in the constraint set $Q$ from an initial state $\boldsymbol{y}_0$, given as
\begin{equation}
    S(c,t,Q,\boldsymbol{y}_0) = \{ \boldsymbol{y}: \boldsymbol{y} \in S(c,t), \gamma\boldsymbol{y} + \boldsymbol{y}_0 \notin Q, 0 \leq \gamma \leq 1 \}. \label{eq: S(c,t,Q,y0)}
\end{equation}
The inclusion of $\boldsymbol{y}_0$ allows $S(c,t,Q,\boldsymbol{y}_0)$ to be generalizable to any starting point in pseudostate-space, not just the origin. The constraint $\gamma\boldsymbol{y} + \boldsymbol{y}_0 \notin Q$ is technically more restrictive than the constraint space $Q$ alone, as it removes any maneuver direction that enters $Q$ at any cost less than or equal to $c$. However, this formulation preserves the linear scaling characteristics of the unconstrained set $S(c,t)$ and can still recover the full array of viable maneuvers on the boundary of $Q$ by varying $c$. For the constraint $\Delta t$, let $T^{\Delta t}(t_{ref})$ be the set of all times in $T$ that are within $\Delta t$ of a given reference time $t_{ref}$, defined as
\begin{equation}
    T^{\Delta t}(t_{ref}) = \{ t: t_{ref} - \Delta t \leq t < t_{ref} + \Delta t, t \in T \}.
\end{equation}
To recover all maneuver times in the control window, consider a sequence of reference times beginning with an arbitrary $t_{ref,0}$ and repeating every $\Delta t$. Let $S(c,T,Q,\boldsymbol{y}_0)$ be the set of all $\boldsymbol{\omega}$ in $S(c,T)$ that are not in the constraint set $Q$ from an initial state $\boldsymbol{y}_0$ for all possible $t_{ref}$, given as
\begin{equation}
    S(c,T,Q,\boldsymbol{y}_0) = \hspace{-4mm} \bigcup_{t_{ref,j} = t_{ref,0} + j\Delta t} \left( \bigcup_{t \in T^{\Delta t}(t_{ref,j})} \hspace{-3.5mm} S(c,t,Q,\boldsymbol{y}_0) \right),
\end{equation}
where $j$ is an integer. This recovers all times in $T$ without repetition for any choice of $t_{ref,0}$. Finally for the constraint $c_{max}$, let $S^*(c,T,c_{max},\Delta t,Q,\boldsymbol{y}_0)$ be the set of $\boldsymbol{\omega}$ that can be reached by $k \geq 1$ control actions of combined cost no greater than $c$ and individual cost no greater than $c_{max}$ at times in $T$ no more frequent than $\Delta t$ that do not enter $Q$ from initial state $\boldsymbol{y}_0$, defined equivalently as
\begin{align}
    S^*(c,T,&c_{max},\Delta t,Q,\boldsymbol{y}_0) \nonumber \\
    = \{ &\boldsymbol{z}: \boldsymbol{z} = \sum\nolimits_{j=1}^{k} \boldsymbol{\Gamma}(t_j)\boldsymbol{u}_j + \boldsymbol{y}_0, \\
    &t_j \in T, | t_{j'} - t_j | \geq \Delta t\hspace{2mm} \forall j' \neq j, \nonumber \\
    &\boldsymbol{u}_j \in U(c_j,t_j), \sum\nolimits_{j=1}^{k} c_j = c, \nonumber \\
    &c_j \leq c_{max}, \gamma\boldsymbol{\Gamma}(t_j)\boldsymbol{u}_j + \boldsymbol{y}_0 \notin Q, 0 \leq \gamma \leq 1 \} \nonumber  \\
    = \{ &\boldsymbol{z}: \boldsymbol{z} = \sum\nolimits_{j=1}^{k} \alpha_j\boldsymbol{y}_j + \boldsymbol{y}_0, \label{eq: S*(c,T,cmax,Q,y0)} \\
    &\boldsymbol{y}_j \in S(c,T^{\Delta t}(t_{ref,j}),Q,\boldsymbol{y}_0), \alpha_j \geq 0, \alpha_j c \leq c_{max}, \nonumber \\
    &t_{ref,j} = t_{ref,0} + j\Delta t_j, \Delta t_j \geq \Delta t, \sum\nolimits_{j=1}^k \alpha_j = \alpha \}. \nonumber
\end{align}
Similar to the unconstrained case in (\ref{eq: S* convex hull}), (\ref{eq: S*(c,T,cmax,Q,y0)}) is defined as the convex hull of $S(c,T,Q,\boldsymbol{y}_0)$ assuming that $\boldsymbol{\omega}$ is reachable. The additional constraints in (\ref{eq: S*(c,T,cmax,Q,y0)}) require that reachability be guaranteed not only by the controllability of the dynamic system in pseudostate space, but also the magnitude constraint $c_{max}$, frequency constraint $\Delta t$, and state constraint $Q$. It is worth noting that this formulation only ensures that one maneuver ($j=1$) will not enter the constraint area $Q$ from the given initial state $\boldsymbol{y}_0$. This issue is addressed by extending the intersecting reachable set approach in (\ref{eq: S* time-ordered}) to (\ref{eq: S*(c,T,cmax,Q,y0)}). Consider subsets of the control window $T$ now defined for a particular time and spacing, given as
\begin{align}
    T^{\Delta t+}_j &= T^{\Delta t+}(t_j) = \{t: t \geq t_j + \Delta t, t \in T \}, \label{eq: Tdt+} \\
    T^{\Delta t-}_j &= T^{\Delta t-}(t_j) = \{t: t \leq t_j - \Delta t, t \in T \}. \label{eq: Tdt-}
\end{align}
Using the reformulation of $S^*(c,T)$ in (\ref{eq: S* time-ordered}), (\ref{eq: S*(c,T,cmax,Q,y0)}) can be defined forward in time as
\begin{align}
    S^*(c,T,c_{max}&,\Delta t,Q) = \{ \boldsymbol{z}: \boldsymbol{z} = \sum\nolimits_{i=1}^{l} \boldsymbol{z}^\Delta_i, \\
    &\boldsymbol{z}^\Delta_i \in S^*(c_i,T^{\Delta t (\cdot)}_i,c_{max},\Delta t,Q,\sum\nolimits_{i'=0}^{i-1} \boldsymbol{z}_{i'}) \nonumber  \\
    &\boldsymbol{z}^\Delta_0 = \Vec{\boldsymbol{0}}, T^{\Delta t (\cdot)}_i \in T^{\Delta t (\cdot)}_{i-1}, \sum\nolimits_{i=1}^l c_i = c \}, \nonumber
\end{align}
where $T^{\Delta t (\cdot)}_i = T^{\Delta t+}_i$ or $T^{\Delta t (\cdot)}_i = T^{\Delta t-}_i$ depending on the chosen direction of time-ordering. Finally, each intersecting reachable set can be limited to a single maneuver without loss of generality, allowing $S^*(c,T,c_{max},\Delta t,Q)$ to be defined as
\begin{align}
    S^*(c,T,c_{max},\Delta t,&Q) = \{ \boldsymbol{z}: \boldsymbol{z} = \sum\nolimits_{i=1}^{l} \boldsymbol{y}_i, \label{eq: S*(c,T,Q)} \\
    &\boldsymbol{y}_i \in S(c_i,T^{\Delta t (\cdot)}_i,Q,\sum\nolimits_{i'=0}^{i-1}\boldsymbol{y}_{i'}) \nonumber  \\
    &\boldsymbol{y}_0 = \Vec{\boldsymbol{0}}, T^{\Delta t (\cdot)}_i \in T^{\Delta t (\cdot)}_{i-1}, \sum\nolimits_{i=1}^l c_i = c \nonumber \\
    & c_i \leq c_{max} \}. \nonumber
\end{align}
This formulation guarantees that constraints $Q$, $\Delta t$, and $c_{max}$ will be satisfied for all maneuvers in the convex hull. It is worth noting here that there is coupling effect between $\Delta t$ and $c_{max}$. Given a $\Delta t = 0$, multiple maneuvers are allowed to exist at a given time $t$, and $c_{max}$ does not restrict the total maneuver magnitude at $t$, only affecting the maneuvers that can enter $Q$ at cost $c_{max}$. The hardware maximum allowable thrust output constraint that motivates the use of the $c_{max}$ can easily be recovered by setting $\Delta t$ to any $\Delta t > 0$.

Returning to the sequential sum of convex hulls in (\ref{eq: RST sequential 1}) for control without operational constraints, the state, time, and magnitude constrained optimal control problem can be expressed through (\ref{eq: S*(c,T,Q)}) as
\begin{align}
    \text{maximize: }& \sum_{i=1}^l \beta_i c \label{eq: RST constrained sequential} \\
    \text{subject to: }& \sum_{i=1}^l \max_{\boldsymbol{y}_i \in S(\beta_ic,T_i^{\Delta t (\cdot)},Q,\sum\nolimits_{i'=0}^{i-1}\boldsymbol{y}_{i'})} 
\hspace{-4mm} \boldsymbol{\lambda}^T\boldsymbol{y}_i \leq \sum_{i=1}^l \beta_i\boldsymbol{\lambda}^T\boldsymbol{\omega}  \nonumber \\
    & \beta_i c = c_i \leq c_{max}, \nonumber
\end{align}
where $T^{\Delta t (\cdot)}_i \in T^{\Delta t (\cdot)}_{i-1}$, $\beta_i \geq 0$, and $\sum_{i=1}^l \beta_i = 1$. Using the same arguments as (\ref{eq: RST SOCP sub-problem}) previously, (\ref{eq: RST constrained sequential}) can be broken into optimal sub-problems, where each sub-problem $i$ is given as
\begin{align}
    \text{maximize: }& \beta_i \boldsymbol{\lambda}^T \boldsymbol{\omega} \label{eq: RST constrained sub-problem} \\
    \text{subject to: }& \max_{t \in T_i^{\Delta t (\cdot)}} g_{S(1,t,Q, \sum_{i'=0}^{i-1} \boldsymbol{y}_{i'})}(\boldsymbol{\lambda}) \leq 1 \nonumber \\
    & \beta_i c = c_i \leq c_{max}. \nonumber
\end{align}
The contact function uses set $S(1,t,c_{max},Q, \sum_{i'=0}^{i-1} \boldsymbol{y}_{i'})$ rather than $U(1,t)$ to account for the location of the state constraint relative both to the maneuver origin $\boldsymbol{y}_0$ and the maximum magnitude allowed $c_{max}$. 

\subsection{Optimality Conditions}

\begin{figure} [b]
    \centering
    \includegraphics[width=1\linewidth]{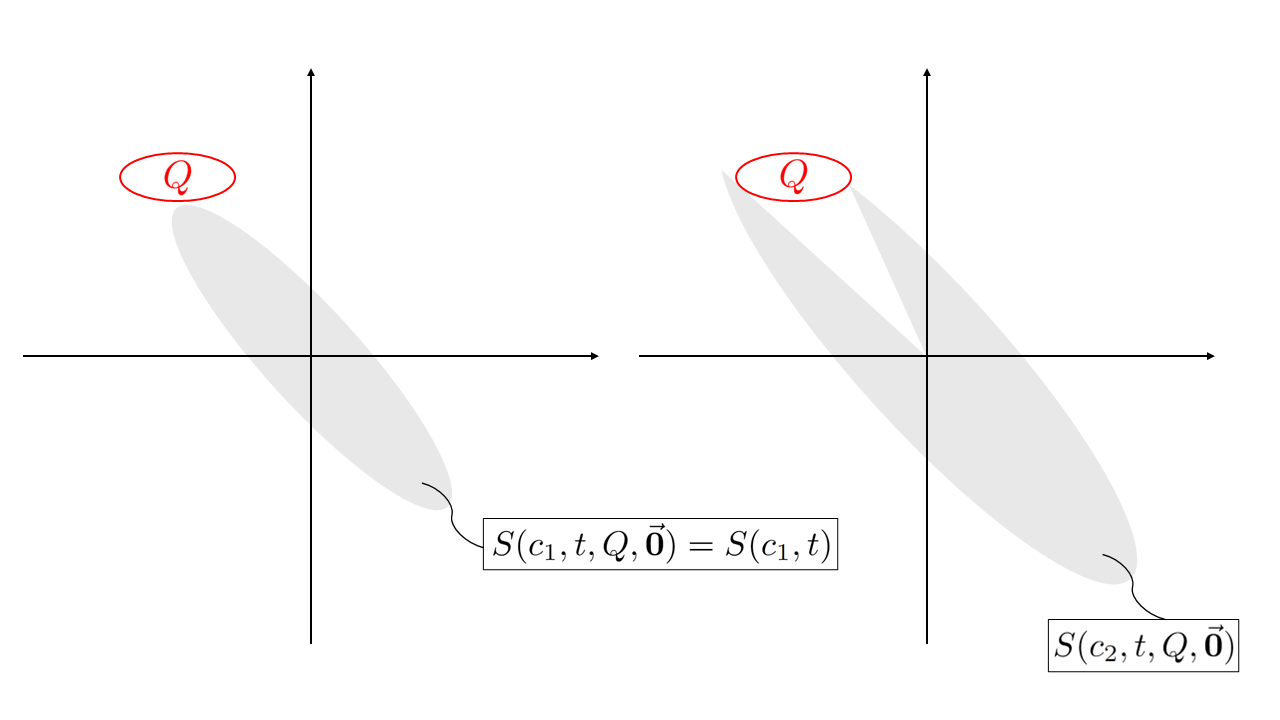}
    \caption{Comparison of the set $S(c,t,Q,\Vec{\boldsymbol{0}})$ for different costs $c_1$ and $c_2$. For a relatively small $c_1$, $Q$ is not reachable for any maneuver direction, and $S(c_1,t,Q,\Vec{\boldsymbol{0}}) = S(c_1,t)$ (left). For a larger $c_2$, $Q$ is reachable for certain maneuver directions, which are removed from $S(c_2,t,Q,\Vec{\boldsymbol{0}})$ (right). }
    \label{fig: constrained S}
\end{figure}

The necessary and sufficient optimality conditions for unconstrained maneuvering in (\ref{eq: g check}), (\ref{eq: efficiency check 1} - \ref{eq: efficiency check 2}), and (\ref{eq: min cost check}) can be generalized to (\ref{eq: RST constrained sub-problem}). However, the inclusion of constraints may result in maneuver directions that are part of $U(1,t)$ but are not allowable states in $S(1,t,Q, \sum_{i'=0}^{i-1} \boldsymbol{y}_{i'})$ for any choice of $\beta_i c$. Consider a relatively small maneuver cost $c_i = \epsilon$. Assuming $c_{max} >> \epsilon$ and $Q$ is not close to $\boldsymbol{y}_0$, the optimal maneuvers found by (\ref{eq: RST constrained sub-problem}) are equivalent to those found in (\ref{eq: RST SOCP sub-problem}). In contrast, assume $c_i = c_{max}$ is infinitely large. The optimal maneuvers found by (\ref{eq: RST constrained sub-problem}) may be significantly different to those found in (\ref{eq: RST SOCP sub-problem}) due to the maneuvers removed in the direction of $Q$. Therefore, the optimal maneuver directions found at a given waypoint are dependent on the potential magnitude of the optimal maneuver. This quality is demonstrated in Figure \ref{fig: constrained S}.

An additional consequence of this maneuver magnitude dependence is the evaluation of $g_{S(1,t,Q, \sum_{i'=0}^{i-1} \boldsymbol{y}_{i'})}(\boldsymbol{\lambda})$ in (\ref{eq: RST constrained sub-problem}). The simple Cauchy-Schwartz evaluation of the contact function in the SOCP constraint in (\ref{eq: RST SOCP}) cannot be applied to (\ref{eq: RST constrained sub-problem}), as maneuvers cannot be guaranteed to be viable without knowing the maneuver magnitude a-priori. Rather than expressing $g_{S(1,t,Q, \sum_{i'=0}^{i-1} \boldsymbol{y}_{i'})}(\boldsymbol{\lambda})$ analytically, the contact function must be evaluated by sampling maneuvers on the surface of the unit spheroid. However, in a similar conclusion to Koenig et al., this sampling approach to solve (\ref{eq: RST constrained sub-problem}) directly is numerically intensive and impractical in application \cite{koenig_RST}. Alternatively, consider expressing the optimality conditions of (\ref{eq: RST constrained sub-problem}) implicitly as a bound on the change in maneuver optimality due to constraints, given as
\begin{align}
    g_{S(c_i,t,Q, \sum_{i'=0}^{i-1}\boldsymbol{y}_{i'})}(\boldsymbol{\lambda})& \label{eq: optimality bound} \\ 
    = g_{U(c_i,t)}&(\boldsymbol{\Gamma}^T(t) \boldsymbol{\lambda}_{opt,i}) - \epsilon_{subopt}, \nonumber
\end{align}
where $\boldsymbol{\lambda}_{opt,i}$ is the solution to the unconstrained sub-problem in (\ref{eq: RST SOCP sub-problem}) and $\epsilon_{subopt}$ bounds the efficiency or optimality difference between constrained and unconstrained optimal maneuvers in sub-problem $i$ of cost $c_i$. The most efficient maneuver possible from a given waypoint $\sum_{i'=0}^{i-1}\boldsymbol{y}_{i'}$ is that found for unconstrained maneuvering, such that the maximum decrease in optimality resulting from constraints can be quantified by sampling valid constrained maneuvers and comparing the resulting states $\boldsymbol{y}$ to the unconstrained $\boldsymbol{\lambda}_{opt,i}$. Because $\epsilon_{subopt}$ is an optimality bound and not an explicit definition of optimality, any arbitrary number of $\boldsymbol{y}$ samples may be used to evaluate $g_{\boldsymbol{y} \in S(1,t,Q, \sum_{i'=0}^{i-1}\boldsymbol{y}_{i'})}(\boldsymbol{\lambda})$. Additional samples may improve the bound at the cost of computational efficiency, introducing a tunable trade-off. $g_{S(1,t,Q, \sum_{i'=0}^{i-1}\boldsymbol{y}_{i'})}(\boldsymbol{\lambda})$ and $\text{arg }g_{S(1,t,Q, \sum_{i'=0}^{i-1}\boldsymbol{y}_{i'})}(\boldsymbol{\lambda})$ then evaluate from the maneuver sample $\Hat{\boldsymbol{u}}$ valid under the given constraints that results in the smallest $\epsilon_{subopt}$. It is also clear that if $\epsilon_{subopt} = 0$, this formulation will result in the unconstrained $g_{U(1,t)}$ and $\text{arg } g_{U(1,t)}$ if $\beta_i c \leq c_{max}$ and the resulting state $\boldsymbol{y}$ that does not enter $Q$.

The necessary and sufficient conditions for sequential constrained optimality can now be derived. Given an arbitrary maneuver direction $\Hat{\boldsymbol{u}}_{i}$ at time $t_{i} \in T_i^{\Delta t (\cdot)}$ at a given magnitude $\beta_i$, the resulting constrained sub-problem $i+1$ is given as
\begin{align}
    \text{maximize: }& \beta_{i+1} \boldsymbol{\lambda}^T \left( \Tilde{\boldsymbol{\omega}} - \beta_i\boldsymbol{\Gamma}(t_i)\Hat{\boldsymbol{u}}_i \right) \label{eq: RST constrained sub-problem i+1} \\
    \text{subject to: }&  \max_{t \in T_{i+1}^{\Delta t(\cdot)}} g_{S(1,t,Q, \sum_{i'=0}^{i-1}\boldsymbol{y}_{i'})}(\boldsymbol{\lambda}) \leq 1. \nonumber
\end{align}
First, $\Hat{\boldsymbol{u}}_{i}$ at time $t_{i}$ is optimal in (\ref{eq: RST constrained sub-problem}) if its contact function satisfies the constraint of (\ref{eq: RST constrained sub-problem i+1}) applying the optimality bound in (\ref{eq: optimality bound}), given as
\begin{equation}
    g_{\Hat{\boldsymbol{u}}_{i}}(\boldsymbol{\Gamma}^T(t_{i})\boldsymbol{\lambda}_{opt,i+1}) \geq 1 - \epsilon_{subopt}, \label{eq: g check constrained}
\end{equation}
for the unconstrained $\boldsymbol{\lambda}_{opt,i+1}$ from (\ref{eq: RST SOCP sub-problem}), such that the contact function for $\Hat{\boldsymbol{u}}_{i}$ is within $\epsilon_{subopt}$ of unconstrained optimality over the subset of times in $T_{i+1}^{\Delta t (\cdot)}$. Second, combining the unconstrained solution $\boldsymbol{\lambda}_{opt,i+1}$ to (\ref{eq: RST constrained sub-problem i+1}) with the optimality bound in (\ref{eq: optimality bound}) yields the optimal maneuver times $t_{opt,i+1}$ and directions $\Hat{\boldsymbol{u}}_{opt,i+1}$, given as
\begin{align}
    t_{opt,i+1} &= \text{arg } \max_{t \in T_{i+1}^{(\cdot)}} g_{S(1,t,Q, \sum_{i'=0}^{i-1}\boldsymbol{y}_{i'})}(\boldsymbol{\lambda}_{opt,i+1}), \\
    \Hat{\boldsymbol{u}}_{opt,i+1} &= \text{arg } g_{S(1,t_{opt,i+1},Q, \sum_{i'=0}^{i-1}\boldsymbol{y}_{i'})}(\boldsymbol{\lambda}_{opt,i+1}), \label{eq: constrained optimal maneuver}
\end{align}
using an arbitrary number of $\boldsymbol{y} \in S(1,t,Q, \sum_{i'=0}^{i-1}\boldsymbol{y}_{i'})$ samples. $\Hat{\boldsymbol{u}}_{i}$ at time $t_{i}$ is optimal in (\ref{eq: RST constrained sub-problem i+1}) if it acts in the same direction and with equal or greater efficiency in pseudostate-space as an optimal maneuver identified in (\ref{eq: constrained optimal maneuver}), defined as
\begin{align}
    &\frac{\left(\boldsymbol{\Gamma}(t_{i})\Hat{\boldsymbol{u}}_{i}\right)^T\left(\boldsymbol{\Gamma}(t_{opt,i+1})\Hat{\boldsymbol{u}}_{opt,i+1}\right)}{||\boldsymbol{\Gamma}(t_{i})\Hat{\boldsymbol{u}}_{i}||_2 ||\boldsymbol{\Gamma}(t_{opt,i+1})\Hat{\boldsymbol{u}}_{opt,i+1}||_2} = 1 \hspace{2mm} \forall \Hat{\boldsymbol{u}}_{opt,i+1}, \label{eq: efficiency check constrained 1} \\
    &\frac{\left(\boldsymbol{\Gamma}(t_{i})\Hat{\boldsymbol{u}}_{i}\right)^T\left(\boldsymbol{\Gamma}(t_{opt,i+1})\Hat{\boldsymbol{u}}_{opt,i+1}\right)}{||\boldsymbol{\Gamma}(t_{opt,i+1})\Hat{\boldsymbol{u}}_{opt,i+1}||^2_2} \geq 1 \hspace{2mm} \forall \Hat{\boldsymbol{u}}_{opt,i+1}, \label{eq: efficiency check constrained 2}
\end{align}
for direction and efficiency respectively. This is the same result as the unconstrained optimality condition in (\ref{eq: efficiency check 1}) and (\ref{eq: efficiency check 2}), with the optimality bound applied implicitly. Third, consider leveraging the definition of cost in the constraint of (\ref{eq: RST SOCP cost}) to realize the remaining cost in a given sub-problem $i$ under the optimality bound in (\ref{eq: optimality bound}), given as
\begin{align}
    \text{maximize: }& \beta_i c \\
    \text{subject to: }& \nonumber \\
    \beta_i c \leq& \frac{\boldsymbol{\lambda}^T \Tilde{\boldsymbol{\omega}}_i}{\max_{t \in T_i^{\Delta t (\cdot)}} g_{S(1,t,Q, \sum_{i'=0}^{i-1} \boldsymbol{y}_{i'})}(\boldsymbol{\lambda})}. \nonumber
\end{align}
$\Hat{\boldsymbol{u}}_{i}$ at time $t_{i}$ is optimal in (\ref{eq: RST constrained sub-problem}) if the minimum cost of sub-problem $i+1$ is no greater than the minimum cost found in the previous sub-problem $i$ minus the cost of the maneuver, equal to $\beta_i$ due to the unitary $\Hat{\boldsymbol{u}}_{i}$, defined as
\begin{align}
    &\frac{\boldsymbol{\lambda}_{opt,i}^T \Tilde{\boldsymbol{\omega}}}{\max_{t \in T_i^{\Delta t (\cdot)}} g_{S(1,t,Q, \sum_{i'=0}^{i-1} \boldsymbol{y}_{i'})}(\boldsymbol{\lambda}_{opt,i})} - \beta_i \label{eq: min cost constrained check} \\
    &\geq \frac{\boldsymbol{\lambda}_{opt,i+1}^T \left( \Tilde{\boldsymbol{\omega}} - \beta_i \boldsymbol{\Gamma}(t_i)\Hat{\boldsymbol{u}}_i \right)}{\max_{t \in T_{i+1}^{\Delta t (\cdot)}} g_{S(1,t,Q, \sum_{i'=0}^{i-1} \boldsymbol{y}_{i'})}(\boldsymbol{\lambda}_{opt,i+1})}. \nonumber
\end{align}

\section{Computationally Efficient Solver for Distributed Space Systems}

The constrained control solver proposed by this work will demonstrate the capabilities of the formulation in (\ref{eq: RST constrained sub-problem}) while preserving computational efficiency for real-time applications. This solver, denoted in this work as the "HD solver," is built upon a solver designed by Koenig et al. for unconstrained control through reachable set theory, denoted in this work as the "KD solver" \cite{koenig_RST}. The KD solver consists of three major components: Initialization, Refinement, and Control Input Extraction. Initialization uses the direction of $\boldsymbol{\omega}$ as a good initial guess of the solution $\boldsymbol{\lambda}_{opt}$ of the SOCP (\ref{eq: RST SOCP}) to produce a set of initial candidate maneuvers times $T_{est}$, defined as the times where (\ref{eq: g_U}) reaches a large value relative to other times in $T$. Refinement iteratively solves (\ref{eq: RST SOCP}) to generate a new $\boldsymbol{\lambda}_{opt}$ for the times in $T_{est}$. This new $\boldsymbol{\lambda}_{opt}$ is used to remove sub-optimal candidate maneuver times from $T_{est}$ whose contact function (\ref{eq: g_U}) falls below a threshold $\epsilon_{remove}$, given as
\begin{equation}
    g_{U(1,t)}(\boldsymbol{\Gamma}^T(t)\boldsymbol{\lambda}_{opt}) < 1 - \epsilon_{remove}, \label{eq: eps remove threshold}
\end{equation}
and add optimal candidate maneuvers times from $T$ to $T_{est}$ that exceed an optimality tolerance $\epsilon_{cost}$, given as
\begin{equation}
    g_{U(1,t)}(\boldsymbol{\Gamma}^T(t)\boldsymbol{\lambda}_{opt}) > 1 + \epsilon_{cost}. \label{eq: eps cost threshold}
\end{equation}
Refinement is repeated until $g_{U(1,t)}(\boldsymbol{\Gamma}^T(t)\boldsymbol{\lambda}_{opt}) < 1 + \epsilon_{cost}$ for all $t \in T_{est}$. Control Input Extraction finds an unconstrained optimal maneuver plan, under only dynamic constraints, at the optimal maneuver times in $T_{est}$ by computing the optimal maneuver directions from (\ref{eq: arg g_U}) and finding the maneuver magnitudes that minimize the resulting pseudostate error of the maneuver plan through a Quadratic Program (QP). Limiting control inputs to times in $T_{est}$ rather than $T$ ensures that maneuvers only occur at optimal times without an explicit definition of optimality in the QP, such that the spacecraft coasts and does not maneuver at sub-optimal times.

The tunable parameters $\epsilon_{remove}$ and $\epsilon_{cost}$ are also incorporated into the HD solver to enhance numerical stability and match the optimality conditions used by the KD solver to generate solutions. Consider the optimality condition derived from the SOCP contact function constraint in (\ref{eq: g check constrained}). This condition matches the structure of (\ref{eq: eps remove threshold}) in the KD solver. Therefore this optimality condition for constrained control sub-problem $i+1$ in the HD solver is given as
\begin{equation}
    g_{\Hat{\boldsymbol{u}}_{i}}(\boldsymbol{\Gamma}^T(t_{i})\boldsymbol{\lambda}_{opt,i+1}) \geq 1 - \epsilon_{remove} - \epsilon_{subopt}, \label{eq: g check constrained alg}
\end{equation}
again considering an arbitrary maneuver direction $\Hat{\boldsymbol{u}}_i$ at time $t_i \in T_i^{\Delta t(\cdot)}$ and a specified $\epsilon_{subopt}$ that bounds the maximum allowable difference between $\Hat{\boldsymbol{u}}_i$ and an optimal unconstrained maneuver. Consider the optimality conditions derived comparatively from the optimal maneuvers $\Hat{\boldsymbol{u}}_{opt,i+1}$ of sub-problem $i+1$ in (\ref{eq: efficiency check constrained 1}) and (\ref{eq: efficiency check constrained 2}). The directional condition should include a $\epsilon_{remove}$ threshold for numerical stability, while the efficiency condition should account for the maximum difference in efficiency between maneuvers deemed optimal under the user-specified tolerance $\epsilon_{cost}$. Therefore the direction and efficiency optimality conditions for constrained control in the HD solver are given as
\begin{align}
    &\frac{\left(\boldsymbol{\Gamma}(t_{i})\Hat{\boldsymbol{u}}_{i}\right)^T\left(\boldsymbol{\Gamma}(t_{opt,i+1})\Hat{\boldsymbol{u}}_{opt,i+1}\right)}{||\boldsymbol{\Gamma}(t_{i})\Hat{\boldsymbol{u}}_{i}||_2 ||\boldsymbol{\Gamma}(t_{opt,i+1})\Hat{\boldsymbol{u}}_{opt,i+1}||_2} \hspace{-.5mm} = \hspace{-.5mm} 1 \hspace{-.5mm} - \hspace{-.5mm} \epsilon_{remove} \hspace{0.5mm} \forall \Hat{\boldsymbol{u}}_{opt,i+1}, \label{eq: efficiency check constrained 1 alg} \\
    &\frac{\left(\boldsymbol{\Gamma}(t_{i})\Hat{\boldsymbol{u}}_{i}\right)^T\left(\boldsymbol{\Gamma}(t_{opt,i+1})\Hat{\boldsymbol{u}}_{opt,i+1}\right)}{||\boldsymbol{\Gamma}(t_{opt,i+1})\Hat{\boldsymbol{u}}_{opt,i+1}||^2_2} \geq \frac{1}{1+\epsilon_{cost}} \hspace{2mm} \forall \Hat{\boldsymbol{u}}_{opt,i+1}. \label{eq: efficiency check constrained 2 alg}
\end{align}
Finally, consider the optimal condition derived from the definition of minimum cost in (\ref{eq: min cost constrained check}). This condition is not directly analogous to a use of $\epsilon_{cost}$ or $\epsilon_{remove}$ in the KD solver. Therefore, the minimum cost optimality condition for constrained control in the HD solver will be the original definition in (\ref{eq: min cost constrained check}).

The HD solver proposed in this work can be understood as two major nested components, Rollout and Branching. Rollout is described by Algorithm \ref{alg: rollout}.
\begin{algorithm} [b]
\caption{Rollout} \label{alg: rollout}
\begin{algorithmic}[1]
    \Statex $\textbf{Input: } t, \Hat{\boldsymbol{u}}, \Tilde{\boldsymbol{\omega}}, T, Q, c_{max}, \Delta c_{min}, \boldsymbol{\epsilon}$
    \State $c_{i} = 0$
    \While{ (\ref{eq: g check constrained alg}), (\ref{eq: efficiency check constrained 1 alg}) and (\ref{eq: efficiency check constrained 2 alg}), or (\ref{eq: min cost constrained check}) and $c_{i} < c_{max}$ }
        \State $c_{i} = c_{i} + \Delta c_{min}$
        \State $\Tilde{\boldsymbol{\omega}} = \Tilde{\boldsymbol{\omega}} - \Delta c_{min}\boldsymbol{\Gamma}(t)\Hat{\boldsymbol{u}}$
        \State $\boldsymbol{\lambda}_{opt,i+1} \gets $ Solve (\ref{eq: RST SOCP sub-problem i+1})
        \State $\Hat{\boldsymbol{u}}_{opt,i+1}(t) \gets $
        \Statex \hspace{15mm} $g_{S(c_{max},T,Q,c_i\boldsymbol{\Gamma}(t)\Hat{\boldsymbol{u}})} (\boldsymbol{\lambda}_{opt,i+1}) \geq 1 - \epsilon_{subopt} $
    \EndWhile
    \Statex $\textbf{Return: } c_{i}, \Hat{\boldsymbol{u}}_{opt,i+1}(t), \boldsymbol{\lambda}_{opt,i+1}$
\end{algorithmic}
\end{algorithm}
Rollout optimizes the magnitude of control input direction $\Hat{\boldsymbol{u}}$ at time $t$ for a desired pseudostate $\Tilde{\boldsymbol{\omega}}$, set of control input times $T$, set of constrained states $Q$, maximum magnitude $c_{max}$, minimum magnitude search step-size $\Delta c_{min}$, and optimality tuning parameters $\boldsymbol{\epsilon} = [\epsilon_{cost}, \epsilon_{remove}, \epsilon_{subopt}]$. Rollout begins by initializing the currently considered sub-problem magnitude $c_{i}$ as zero (Alg. \ref{alg: rollout}, Line 1) and then increases $c_{i}$ to find the remaining $\Tilde{\boldsymbol{\omega}}$ relative to the new waypoint (Alg. \ref{alg: rollout}, Lines 3-4). This waypoint is used to solve the unconstrained control problem in (\ref{eq: RST SOCP sub-problem i+1}) using the KD algorithm, producing $\boldsymbol{\lambda}_{opt,i+1}$ (Alg. \ref{alg: rollout}, Line 5). This unconstrained result is leveraged to bound constrained maneuver optimality, and optimal maneuvers $\Hat{\boldsymbol{u}}_{opt,i+1}(t)$ for the sub-problem $i+1$ are found as maneuvers that do not enter $Q$ at cost $c_{max}$ that are within the sub-optimality bound $\epsilon_{subopt}$ of unconstrained optimality (Alg. \ref{alg: rollout}, Line 6). As mentioned previously, the contact function for constrained maneuvering cannot be evaluated as a single expression through the Cauchy-Schwartz inequality due to the effect of $Q$. Therefore, Rollout evaluates a coarse set of unit sphere samples in $\Re^3$ space, introducing a tradeoff between potential optimality gains found by finer sampling and the increased computational effort required to evaluate the contact function of additional samples. While the coarseness of the unit sphere samples is a hyperparameter of the solver, it does not affect the constrained sub-optimality guarantee $\epsilon_{subopt}$ set by the user, so long as $\Tilde{\boldsymbol{\omega}}$ remains controllable. The algorithm continues to increase $c_{i}$ until none of the optimality conditions in (\ref{eq: g check constrained alg}), (\ref{eq: efficiency check constrained 1 alg}) and (\ref{eq: efficiency check constrained 2 alg}), or (\ref{eq: min cost constrained check}) are true, indicating that the maximum optimal maneuver magnitude has been found (Alg. \ref{alg: rollout}, Line 2). Rollout returns the optimal control input magnitude $c_{i}$ for control input direction $\Hat{\boldsymbol{u}}$ at time $t$, maneuvers $\Hat{\boldsymbol{u}}_{opt,i+1}(t)$ viable under constraint $Q$ and within the constrained sub-optimality bound $\epsilon_{subopt}$, and the outward normal $\boldsymbol{\lambda}_{opt,i+1}$ of $S^*(c,T)$ of the resulting unconstrained sub-problem $i+1$. In terms of computational load, Rollout requires at most $c_{max}/\Delta c_{min} + 1$ iterations (rounded up to an integer value), requiring $c_{max}/\Delta c_{min} + 1$ unconstrained optimal solutions from the KD solver (Alg. \ref{alg: rollout}, Line 5). Rollout is the primary building block used to build Branching, described by Algorithm \ref{alg: branching}.

\begin{figure}[b]
    \centering
    \includegraphics[width=1\linewidth]{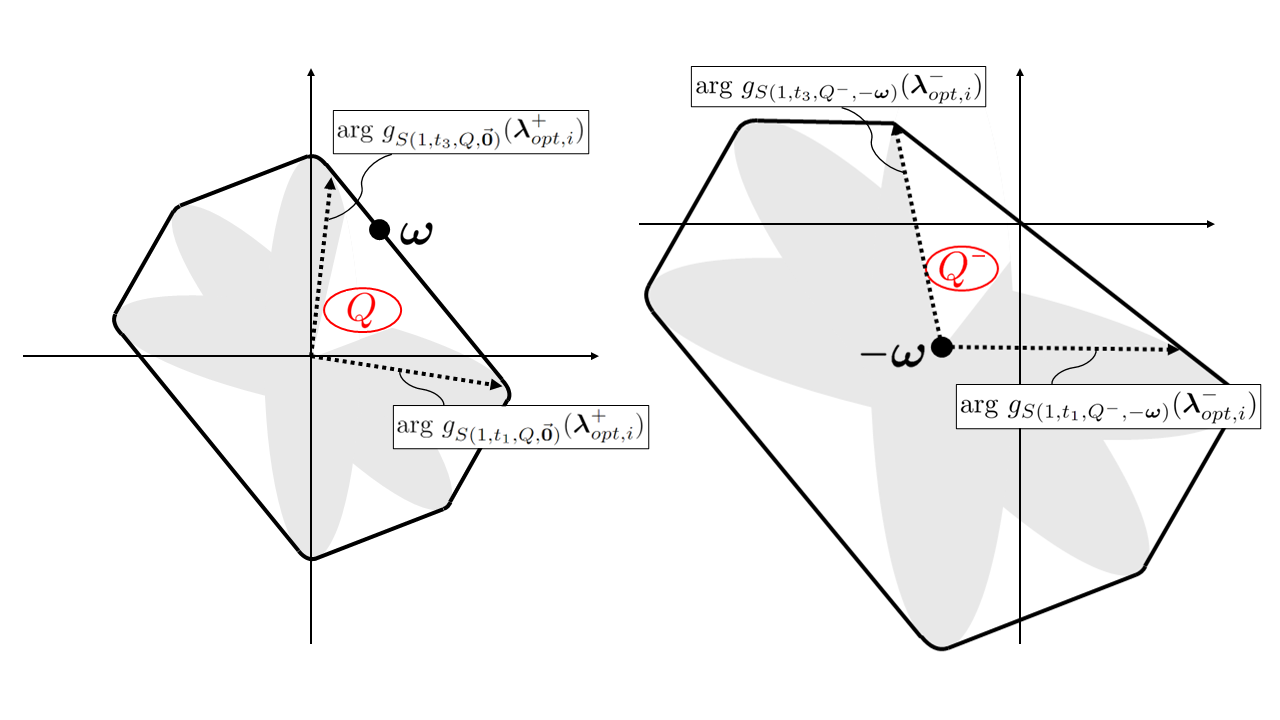}
    \caption{Reachable sets and optimal maneuver directions under a symmetric cost function constrained by $Q$ from a waypoint at the origin (left) and a waypoint at the desired pseudostate (right). Reachable set differences are caused by the relative location of $Q$ to the considered waypoint. }
    \label{fig: forward and backward}
\end{figure}

\begin{algorithm} [t]
\caption{Branching} \label{alg: branching}
\begin{algorithmic} [1]
    \Statex $\textbf{Input: } \boldsymbol{\omega}, T, Q, c_{max}, \Delta c_{min}, \Delta t, \boldsymbol{\epsilon}$\
    \State Initialize $\boldsymbol{u}(t)$ empty, $c_{subopt} = 0$, $\Tilde{\boldsymbol{\omega}} = \boldsymbol{\omega}$
    \State $\boldsymbol{u}_{uncon}(t) \leftarrow$ Solve unconstrained SOCP (\ref{eq: RST SOCP})
    \While{ $\boldsymbol{u}_{uncon}(t)$ enters $Q$, exceeds $c_{max}$, or is more frequent than $\Delta t$}
        \State $\Hat{\boldsymbol{u}}_{i}^+(t), \boldsymbol{\lambda}^+_{opt,i} \gets$ Rollout($0,\Vec{\boldsymbol{0}}, \Tilde{\boldsymbol{\omega}}, T, Q, ...)$
        \State $t^+_i, \Hat{\boldsymbol{u}}^+_i \gets \Hat{\boldsymbol{u}}^+(t) \text{ closest to } t_0$
        \State $c^+_i, \Hat{\boldsymbol{u}}^+(t), \boldsymbol{\lambda}^+_{opt,i+1} \hspace{-1.5mm} \gets$ \hspace{-1mm} Rollout($t^+_i, \Hat{\boldsymbol{u}}^+_i, \Tilde{\boldsymbol{\omega}}, T, Q, ...$)
        \State Add $t^+_i, c^+_i\Hat{\boldsymbol{u}}^+_i$ to $\boldsymbol{u}(t)$
        \State $c_{subopt} = c_{subopt} + \boldsymbol{\lambda}^+_{opt,i+1}(\Tilde{\boldsymbol{\omega}} - c^+_i\boldsymbol{\Gamma}(t^+_i)\Hat{\boldsymbol{u}}^+_i) $
        \Statex \hspace{25mm}$  - (\boldsymbol{\lambda}^+_{opt,i}\Tilde{\boldsymbol{\omega}} - c^+_i) + \Delta c_{min} $
        \State $\Tilde{\boldsymbol{\omega}} = \Tilde{\boldsymbol{\omega}} - c^+_i\boldsymbol{\Gamma}(t^+_i)\Hat{\boldsymbol{u}}^+_i, T \gets T^{\Delta t+}_i$
        \State $Q = \{\boldsymbol{x}': \boldsymbol{x}' = \boldsymbol{x} - c^+_i\boldsymbol{\Gamma}(t^+_i)\Hat{\boldsymbol{u}}^+_i, \boldsymbol{x} \in Q \}$
        
        \State $Q^- = \{\boldsymbol{x}': \boldsymbol{x}' = \boldsymbol{x} - \Tilde{\boldsymbol{\omega}}, \boldsymbol{x} \in Q \}$
        \State $\Hat{\boldsymbol{u}}_{i}^-(t), \boldsymbol{\lambda}^-_{opt,i} \gets$ Rollout($0,\Vec{\boldsymbol{0}}, \text{-}\Tilde{\boldsymbol{\omega}}, T, Q^-, ...)$
        \State $t^-_i, \Hat{\boldsymbol{u}}^-_i \gets \Hat{\boldsymbol{u}}^-(t) \text{ closest to } t_f$
        \State $c^-_i, \Hat{\boldsymbol{u}}^-(t), \boldsymbol{\lambda}^-_{opt,i+1} \hspace{-1.5mm} \gets$ \hspace{-1.5mm} Rollout($t^-_i, \Hat{\boldsymbol{u}}^-_i, \text{-}\Tilde{\boldsymbol{\omega}}, T, Q^-, ...$)
        \State Add $t^-_i, c^-_i\Hat{\boldsymbol{u}}^-_i$ to $\boldsymbol{u}(t)$
        \State $c_{subopt} = c_{subopt} + \boldsymbol{\lambda}^-_{opt,i+1}(\Tilde{\boldsymbol{\omega}} - c^-_i\boldsymbol{\Gamma}(t^-_i)\Hat{\boldsymbol{u}}^-_i) $
        \Statex \hspace{25mm}$ - (\boldsymbol{\lambda}^-_{opt,i}\Tilde{\boldsymbol{\omega}} - c^-_i) + \Delta c_{min} $
        
        \State $\Tilde{\boldsymbol{\omega}} = \Tilde{\boldsymbol{\omega}} - c^-_i\boldsymbol{\Gamma}(t^-_i)\Hat{\boldsymbol{u}}^-_i, T \gets T^{\Delta t-}_i$
        \State $\boldsymbol{u}_{uncon}(t) \leftarrow$ Solve unconstrained SOCP (\ref{eq: RST SOCP})
    \EndWhile
    \State Add $\boldsymbol{u}_{uncon}(t)$ to $\boldsymbol{u}(t)$
    \Statex $\textbf{Return: } \boldsymbol{u}(t), c_{subopt}$
\end{algorithmic}
\end{algorithm}
Branching uses Rollout to perform sequential constrained maneuvering until the remaining control problem described by $\Tilde{\boldsymbol{\omega}}$ and $T$ can be solved with unconstrained maneuvering. Here, the inputs of Rollout are abbreviated for clarity to only show those that vary between each iteration of the algorithm. Branching begins by initializing an empty maneuver plan $\boldsymbol{u}(t)$, residual sub-optimality bound $c_{subopt}$ as zero, and waypoint pseudostate $\Tilde{\boldsymbol{\omega}}$ as the initial pseudostate $\boldsymbol{\omega}$ and finds a valid maneuver plan for unconstrained maneuvering using the KD solver (Alg. \ref{alg: branching}, Lines 1-2). The residual sub-optimality bound $c_{subopt}$ is introduced to quantify the maximum possible cost difference between the realized constrained maneuver plan and a globally-optimal solution. Branching finds constrained optimal maneuvers sequentially forward and backward in time until the unconstrained solution is valid under the constraints $Q, c_{max}, $ and $\Delta t$ (Alg. \ref{alg: branching}, Lines 3-19). For forward time-ordered sequential maneuvering, Branching uses Rollout with a dummy zero maneuver to find the maneuver times and directions valid under $Q$ and optimal given $\epsilon_{subopt}$, and the corresponding unconstrained optimal $\boldsymbol{\lambda}^+_{opt,i}$ for a waypoint at the origin (Alg. \ref{alg: branching}, Line 4). The algorithm then identifies a single forward time-ordered branch as the earliest constrained optimal maneuver (Alg. \ref{alg: branching}, Line 5). The HD solver proposed by this work only considers a single branch at each waypoint, rather than an exhaustive search of all optimal maneuver branches at each waypoint, to conserve computational efficiency. Using the maneuver that constrains the fewest future maneuver times limits the impact of this strategy on optimality or controllability. Rollout is then run for this maneuver to find the optimal maneuver magnitude, the optimal maneuvers of the next sub-problem $i+1$ valid under $Q$ at the waypoint resulting from this maneuver magnitude, and the unconstrained optimal $\boldsymbol{\lambda}^+_{opt,i+1}$ (Alg. \ref{alg: branching}, Line 6). The maneuver at its optimal magnitude is added to the maneuver plan, and the residual sub-optimality of the maneuver is added to the overall control problem bound (Alg. \ref{alg: branching}, Line 7-8). This sub-optimality bound is taken as the difference in the remaining unconstrained cost between sub-problems $i$ and $i+1$, accounting for the maneuver cost $c^+_i$, and added to the user-defined minimum maneuver magnitude step-size $\Delta c_{min}$ used in Rollout. The remaining pseudostate and candidate maneuver times are updated, and $Q$ is adjusted to be defined relative to this new waypoint (Alg. \ref{alg: branching}, Lines 9-10). This process is repeated similarly for backward time-ordered maneuvering by starting from the waypoint at $\Tilde{\boldsymbol{\omega}}$, adjusting the relative location of $Q$ to $Q^-$, using Rollout with a dummy zero maneuver to find the maneuver times and directions valid under $Q^-$, and taking the latest constrained maneuver as the single branch (Alg. \ref{alg: branching}, Line 11-13). Rollout is performed on this branch, the optimized maneuver is added to the maneuver plan, and the residual sub-optimality and remaining pseudostate and candidate maneuver times are all updated (Alg. \ref{alg: branching}, Line 14-17). It is worth noting that the backward time-ordered problem here poses the control problem as the reverse of the forward time-ordered problem, using $-\Tilde{\boldsymbol{\omega}}$ and adjusting $Q^-$ relative to $-\Tilde{\boldsymbol{\omega}}$. This reverse setup is valid only if the cost function is not only norm-like but also symmetric, such that a maneuver is equally as efficient in both the parallel and anti-parallel directions. However, most norm-like cost functions used in spacecraft applications are symmetric, and this assumption is made without significant loss of applicability. The relationship between the constrained forward and backward time-ordered problems is visualized in Figure \ref{fig: forward and backward}.
An optimal unconstrained maneuver plan is found using the KD algorithm for the reduced control problem after applying both maneuvers, and Branching iterates until the unconstrained maneuver plan satisfies the constraints (Alg. \ref{alg: branching}, Line 18). The unconstrained maneuver plan is then added to the output maneuver plan (Alg. \ref{alg: branching}, Line 20). Branching returns the constrained maneuver plan $\boldsymbol{u}(t)$ and the optimality bound $c_{subopt}$ between this maneuver plan and the actual optimal cost of the problem. In terms of computational load, each iteration of Branching requires a solution of the unconstrained problem at each waypoint (Alg. \ref{alg: branching}, Lines 4 and 12), an optimization of each forward and backward time-ordered constrained maneuver magnitude using Rollout (Alg. \ref{alg: branching}, Lines 6 and 14), and a single unconstrained maneuver plan (Alg. \ref{alg: branching}, Line 18) to add two constrained maneuvers to the maneuver plan. If $n$ constrained maneuvers are required to produce a valid maneuver plan, Branching requires at most $n(c_{max}/\Delta c_{min} + 1) + 3n/2 + 1$ unconstrained optimal solutions from the KD solver.

\section{Validation}

The sequential reachable set theory formulation and proposed HD solver for constrained impulsive control are validated on general DSS use-case reconfigurations. DSS reconfigurations consider the relative motion of a spacecraft, known as the deputy and denoted with subscript $d$, linearized about a reference orbit of a virtual or real non-maneuvering spacecraft, known as the chief and denoted as subscript $c$. This work implements an LTV relative motion model using two state represenations: Radial-Tangential-Normal (RTN) and quasi-nonsingular Relative Orbital Elements (ROE) \cite{damico_phd}. RTN is a cartesian frame whose origin is defined at the chief and whose axes are aligned radially from the central body of gravitation, normally along the chief's angular momentum vector, and tangentially to complete the right-handed triad.  RTN is commonly used to describe deputy control actions. ROE is an Integration Constant state representation valid for orbits of arbitrary eccentricity, formulated as the nonlinear combination of Keplerian orbital elements and defined as
\begin{equation} \label{eq: ROE}
    \boldsymbol{x} = 
    \begin{bmatrix}
        \delta a \\ \delta \lambda \\ \delta \boldsymbol{e} \\ \delta \boldsymbol{i}
    \end{bmatrix} = 
    \begin{bmatrix}
    \delta a \\ \delta \lambda \\ \delta e_x \\ \delta e_y \\ \delta i_x \\ \delta i_y
    \end{bmatrix} = 
    \begin{bmatrix}
    \Delta a/a_c \\
    \Delta u + \Delta \Omega \cos i_c \\
    e_d\cos\omega_d - e_c\cos\omega_c \\
    e_d\sin\omega_d - e_c\sin\omega_c \\
    \Delta i \\
    \Delta \Omega \sin i_c
    \end{bmatrix},
\end{equation} 
\noindent where $\Delta$ is a difference in the associated quantity between the deputy (subscript $d$) and the chief (subscript $c$), $\delta a$ is the relative semi-major axis, $\delta \lambda$ is the relative mean longitude, $\delta e_x$ and $\delta e_y$ are the components of the relative eccentricity vector $\delta \boldsymbol{e}$, $\delta i_x$ and $\delta i_y$ are the components of the relative inclination vector $\delta \boldsymbol{i}$, $u = M + \omega$ is the mean argument of latitude, and $a$, $e$, $i$, $\Omega$, $\omega$, and $M$ are the Keplerian orbital elements. The Validation scenarios all use mean ROE computed from mean Keplerian orbital elements, where "mean" refers to the orbit average of J2-perturbed two-body motion to remove "osculating" short-period oscillations. The full LTV relative motion model of DSS is described in Appendix A.

Each of the following Validation scenarios will share several base parameters of the reconfiguration problem. The initial absolute orbit of the chief spacecraft at $t_0$ in all scenarios is given in Table \ref{tab: chief oe0}.
\begin{table}[b]
    \centering
    \caption{\bf Initial Mean Quasi-Nonsingular Chief Orbit}
    \begin{tabular}{c c c c c}
        $a(km)$ & $e_x(-)$ & $e_y(-)$ & $i(^\circ)$ & $\Omega(^\circ)$ \\
        \hline \hline 
        \rule{0pt}{2ex} 25000 & -0.658 & -0.239 & 51 & 30 \\
        \hline
    \end{tabular}
    \label{tab: chief oe0}
\end{table}
The initial mean argument of latitude $u_0$ is varied between each scenario to produce clear pseudostate trajectories of differing reconfigurations and constraints for visual analysis. All reconfiguration scenarios will use a 39000 second control window, approximately one orbital period of the chief orbit. The control domain $T$ is taken uniformly every ten seconds relative to $t_0$, such that
\begin{equation}
    T = \{ 0, 10, 20, ... 39000\}.
\end{equation}
This results in 3901 candidate maneuver times. The chosen cost function is the sum of L2-norms of RTN control actions at times in $T$, and control is admissible at all times in the control window. This simple cost function, as opposed to a time-varying cost function, does not impact the generality of the definition in (\ref{eq: cost function}) and is used to allow the specific advantages of the formulation proposed in this work to be more readily apparent.

Validation is broken into three sections, Sequential Formulation, Convex Constraints, and Nonconvex Constraints, to fully demonstrate the capabilities of the approach proposed by this work. Specifically, three impulsive control solvers from the literature will be used comparatively to evaluate the optimality, constraint satisfaction, and computational efficiency of the HD solver. Each algorithm's runtime will be measured on the same hardware, using a 3.20 GHz processor and 16GB RAM. First, the KD solver will generate unconstrained solutions to each Validation scenario using the same optimality tolerances as the HD solver ($\epsilon_{cost} = \epsilon_{remove} = 0.01)$. The Initialization component of the KD solver is set to evaluate the contact function for every ten candidate times in $T$, taking the ten highest as the initial $T_{est}$, and the Control Input Extraction component uses an identity matrix to minimize the pseudostate error. For the Sequential Formulation scenario, the KD solver provides a provably-optimal solution to the original unconstrained reachable set theory problem with which solutions from the HD solver are compared to prove the validity of the sequential approach. The KD solver is also used with the Convex and Nonconvex Constraint scenarios, neglecting the chosen constraints, to provide a visual trajectory comparison between optimal unconstrained and constrained control. Second, the SDPT3 interior point solver within the CVX software package, denoted in this work as the "CVX solver," generates optimal solutions to impulsive control problems under convex constraints, including SOCPs, QPs, Linear Programs (LPs), and semidefinite programs (SDPs)  \cite{boyd}. For the Convex Constraints scenario, the CVX solver will produce an optimal solution to an equivalent primal version of the dual in (\ref{eq: RST constrained sequential}). This primal is given as
\begin{align}
    \text{minimize: }& \sum_{i=1}^l||\boldsymbol{u}_i||_2 \label{eq: cvx problem} \\
    \text{subject to: }& \boldsymbol{z}_i = \boldsymbol{\Gamma}(t_i)\boldsymbol{u}_i + \boldsymbol{z}_{i-1} \nonumber \\
    & t_i \geq t_{i-1}, \boldsymbol{z}_0 = \Vec{\boldsymbol{0}}, \boldsymbol{z}_l = \boldsymbol{\omega} \nonumber\\
    & \boldsymbol{z}_i \notin Q \nonumber\\
    & ||\boldsymbol{u}_i||_2 \leq c_{max}. \nonumber
\end{align}
The maneuver plans generated by the HD solver will be compared to this solution to verify bounded optimality and evaluate computational effort. Third, an SCP optimization algorithm from the literature, denoted in this work as the "SCP solver," is leveraged to find a locally optimal solution to the impulsive control problem under nonconvex constraints \cite{guffanti_ART}. This SCP solver considers a nonconvex Keep Out Ellipsoid (KOE) constraint space, linearized at each iteration $k$ about the nominal trajectory of the last iteration $k-1$ to series of waypoint halfspaces, given as
\begin{align}
    \boldsymbol{a}_{KOE,i}^T &= (\Bar{\boldsymbol{z}}_i - \boldsymbol{x}_c)^T\boldsymbol{D}, \\
    b_{KOE,i} &= \sqrt{(\Bar{\boldsymbol{z}}_i - \boldsymbol{x}_c)^T \boldsymbol{D} (\Bar{\boldsymbol{z}}_i - \boldsymbol{x}_c)},
\end{align}
where $\boldsymbol{a}_{KOE,i}$ and $b_{KOE,i}$ define the linearized halfspace for waypoint $i$, $\Bar{\boldsymbol{z}}_i$ defines a waypoint of the trajectory of the last iteration, $\boldsymbol{x}_c$ is the centroid of the KOE, and $\boldsymbol{D} = \text{diag}(1/r_1^2, 1/r_2^2, ... 1/r_n^2)$ sizes the KOE along each ROE dimension. The SCP solver is extended to incorporate time constraints through a series of maneuver efficiency coefficients $\boldsymbol{\rho}_{\Delta t}$, one coefficient bounded between 0 and 1 for each time in $T$. $\boldsymbol{\rho}_{\Delta t}$ is initialized as a series of 1s and updated each iteration $k$ by subtracting 0.1 from all coefficients within $\Delta t$ of a maneuver in the last iteration's maneuver plan and adding 0.1 to all other coefficients. The difference between the trajectory of the current iteration and the last iteration is bounded by a trust region constraint at each waypoint, given as
\begin{equation}
    ||\boldsymbol{z}_i - \Bar{\boldsymbol{z}}_i||_2 \leq \zeta_k,
\end{equation}
where $\zeta_k$ is the trust region at iteration $k$. The approximate convex problem solved at each iteration of the SCP solver is given as
\begin{align}
    \text{minimize: }& \sum_{i=1}^l||\boldsymbol{u}_i||_2 \label{eq: scp convex approx} \\
    \text{subject to: }& \boldsymbol{z}_i = \rho_{\Delta t,i}\boldsymbol{\Gamma}(t_i)\boldsymbol{u}_i + \boldsymbol{z}_{i-1} \nonumber \\
    & t_i \geq t_{i-1}, \boldsymbol{z}_0 = \Vec{\boldsymbol{0}}, \boldsymbol{z}_l = \boldsymbol{\omega} \nonumber\\
    & \boldsymbol{a}_{KOE,i}^T (\boldsymbol{z}_i - \boldsymbol{x}_c ) \geq b_{KOE,i} \nonumber\\
    & ||\boldsymbol{u}_i||_2 \leq c_{max} \nonumber \\
    & ||\boldsymbol{z}_i - \Bar{\boldsymbol{z}}_i||_2 \leq \zeta_k. \nonumber
\end{align}
The trust region is shrunk each iteration to encourage solver convergence with an exponential shrink factor, given as
\begin{equation}
    \zeta_{k+1} = \left( \frac{\zeta_K}{\zeta_1} \right)^{1/K} \zeta_k,
\end{equation}
where $K$ is the maximum number of iterations and $\zeta_1$ and $\zeta_K$ are the user specified trust regions for iteration 1 and $K$ respectively. The maneuver plan at each iteration is found by taking the largest maneuver from the maneuver plan found by the approximate convex problem, removing all maneuvers at times within $\Delta t$ of this maneuver, and repeating until no maneuver magnitudes are larger than a user-specified tolerance $c_{tol}$. The SCP solver terminates after $K$ iterations, or when $\zeta_k \leq \zeta_K$ and the decrease in cost $c$ between iteration $k$ and $k-1$ is less than $c_{tol}$. For the Nonconvex Constraints scenario, the cost, constraint satisfaction, and required computational effort of the solution found by the SCP solver will be compared against the solution found by the HD solver. 

\subsection{Sequential Formulation}
\begin{figure*}[h!]
    \centering
    \includegraphics[width=1\linewidth]{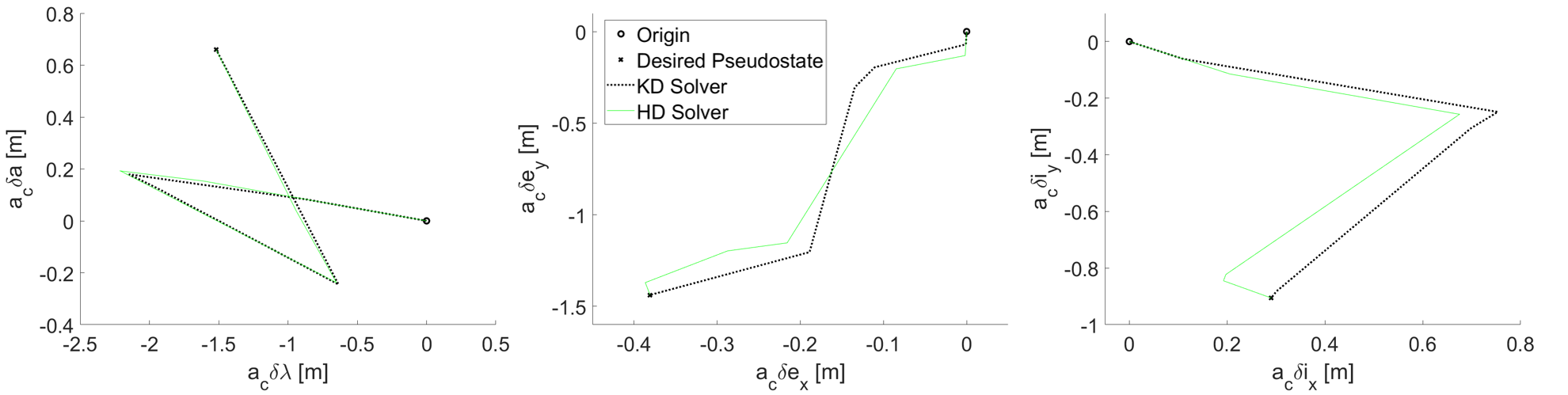}
    \caption{Pseudostate trajectory of the forward time-ordered setup (Alg. \ref{alg: branching}, Lines 4-10) compared to the solution of the KD solver for unconstrained control.}
    \label{fig: sequential formaultion forward}
\end{figure*}
\begin{figure*}[h]
    \centering
    \includegraphics[width=1\linewidth]{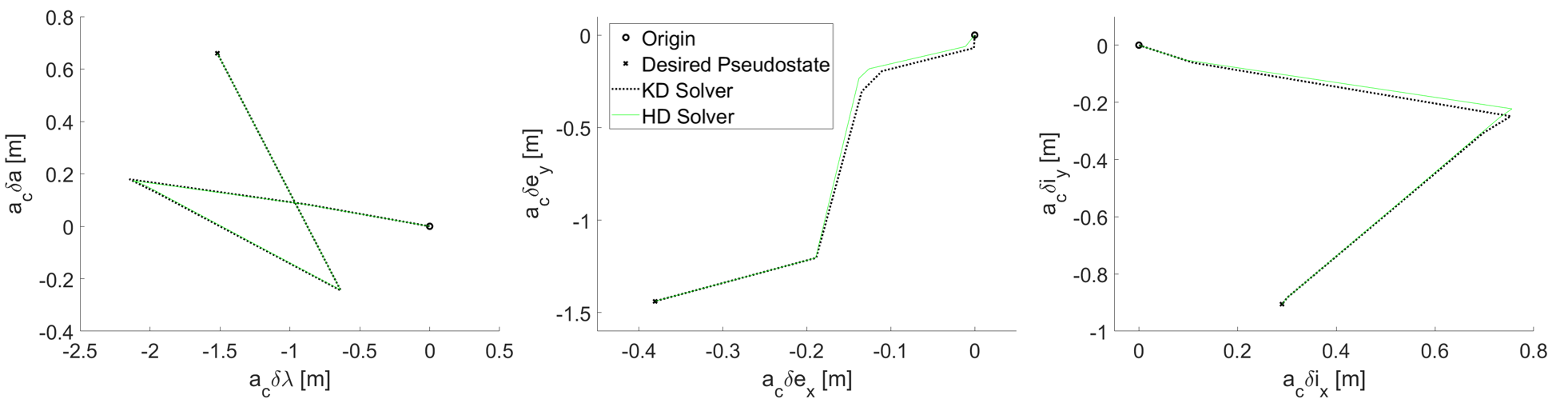}
    \caption{Pseudostate trajectory of the backward time-ordered setup (Alg. \ref{alg: branching}, Lines 11-17) compared to the solution of the KD solver for unconstrained control.}
    \label{fig: sequential formaultion backward}
\end{figure*}

First, the sequential time-ordered formulation and optimal sub-problem structure are validated against the KD solver. The initial and final desired deputy states, and resulting pseudostate, that define this scenario are provided in Table \ref{tab: sequential formulation roe} for an initial chief spacecraft $u = 65^\circ$ at $t_0$.
\begin{table}[b]
    \centering
    \setlength{\tabcolsep}{5.5pt}
    \caption{\bf Initial, Desired, and Pseudostate Mean ROE Deputy States for Sequential Formulation Validation (m)}
    \begin{tabular}{l|c c c c c c}
         & $a_c\delta a$ & $a_c \delta \lambda$ & $a_c \delta e_x$ & $a_c \delta e_y$ & $a_c \delta i_x$ & $a_c \delta i_y$ \\
         \hline \hline
        \rule{0pt}{2ex}$\boldsymbol{x}_0$ & -2.96 & -25.14 & 5.16 & 384.13 & -4.91 & -442.85 \\
        $\boldsymbol{x}_f$ & -2.30 & 0.81 & 4.09 & 383.28 & -4.62 & -442.85 \\
        $\boldsymbol{\omega}$ & 0.66 & -1.52 & -0.38 & -1.44 & 0.29 & -0.91 \\
        \hline
    \end{tabular}
    \label{tab: sequential formulation roe}
\end{table}
This problem is very similar to that posed by Koenig et al. to validate the KD solver \cite{koenig_RST}. Unconstrained control is realized within the HD solver by setting $Q$ as an empty set, $\Delta t$ as 0s, and $c_{max}$ as 1e-3m/s. The minimum step-size $\Delta c_{min}$ of the optimal maneuver magnitude search in Rollout is set as 1e-7m/s. To validate the sequential formulation, the optimal trajectory is generated using forward time-ordering and backward time-ordering in isolation. The forward time-ordered setup neglects Lines 11-18 in Alg. \ref{alg: branching}, while the backward time-ordered setup neglects Lines 4-10 and Line 18 in Alg. \ref{alg: branching}. Both setups modify the loop termination conditions from constraint satisfaction to $\Tilde{\boldsymbol{\omega}} \approx 0$ (Alg. \ref{alg: branching}, Line 3). Because the system is unconstrained, Rollout can evaluate maneuver optimality simply through the Cauchy-Schartz inequality in (\ref{eq: g_U}), and maneuver direction sampling is not required (Alg. \ref{alg: rollout}, Line 6). The resulting ROE trajectory of the HD solver and the KD solver are displayed in Figures \ref{fig: sequential formaultion forward} and \ref{fig: sequential formaultion backward}.

Both the forward and backward time-ordered setups follow a similar, but slightly different, pseudostate trajectory to the KD solver. This is a reasonable result, recalling that an optimal control input profile may not be unique if $\boldsymbol{\omega}$ lies on an edge or point of the boundary of $S^*(c_{opt},T)$. The KD solver finds an unconstrained impulsive 3-maneuver plan that satisfies the tolerance specified by the user through $\epsilon_{cost}$ in 4 iterations with a total $\Delta$v cost of 23.03e-5m/s. The optimal lower bound on the reconfiguration is found through $\boldsymbol{\lambda}^T_{opt}\boldsymbol{\omega}/\max_{t \in T} g_{U(1,t)}(\boldsymbol{\Gamma}^T(t)\boldsymbol{\lambda}_{opt})$ as 22.94e-5m/s. The maneuver plan achieves the desired state with a root mean square error of 23.56e-5m, $<0.01\%$ error normalized to the original pseudostate ($||\boldsymbol{\omega}_{err}||_2/||\boldsymbol{\omega}||_2$). The forward time-ordered setup of the HD solver finds a 10-maneuver plan with a total $\Delta$v cost of 23.08e-5m/s, while the backward time-ordered setup finds a 12-maneuver plan with a total $\Delta$v cost of 23.13e-5m/s. The sub-optimality bounds $c_{subopt}$ for both approaches are found simply as the sum of minimum step-sizes used to optimize the maneuver magnitudes, 10e-7m/s and 12e-7m/s respectively. In this case where $\epsilon_{subopt} = 0$, $c_{subopt}$ represents the worst case increase in cost between the KD solver and the HD solver due to using the sequential approach. Both costs found by the HD solver are within their associated $c_{subopt}$ of the maneuver plan found by the KD solver. Furthermore, the forward time-ordered plan achieves the desired state with a root mean square error of 67.65e-5m, while the backward time-ordered plan has a root mean square error of 86.98e-5m. Both errors are $<0.01\%$ normalized to the original pseudostate. The sequential approach does increase computational effort. The KD algorithm required 1.90s to solve trajectory optimization, while the forward and backward time-ordered setups had 186.06s and 135.91s runtimes. Given the KD solver computational time, a rough upper bound on the expected computation time of 10 and 12 instances of Rollout, one for each maneuver, can be found as 10 maneuvers * 1.90s runtime * $(c_{max}/\Delta c_{min} + 1) \approx$ 190,019s and 12 maneuvers * 1.90s runtime * $(c_{max}/\Delta c_{min} + 1) \approx$ 228,023s. The forward and backward time-ordered setup runtimes are well within these bounds. Therefore, the sequential time-ordered reachable set theory formulation in (\ref{eq: RST sequential 2}) and resulting sub-problems in (\ref{eq: RST sequential sub-problem}) are empirically demonstrated to be equivalent to the original unconstrained reachable set definition in (\ref{eq: S* convex hull}). This conclusion for a single example scenario will be further verified through Monte Carlo simulations of the HD solver for convex and nonconvex constraints.

\begin{figure*}[h]
    \centering
    \includegraphics[width=1\linewidth]{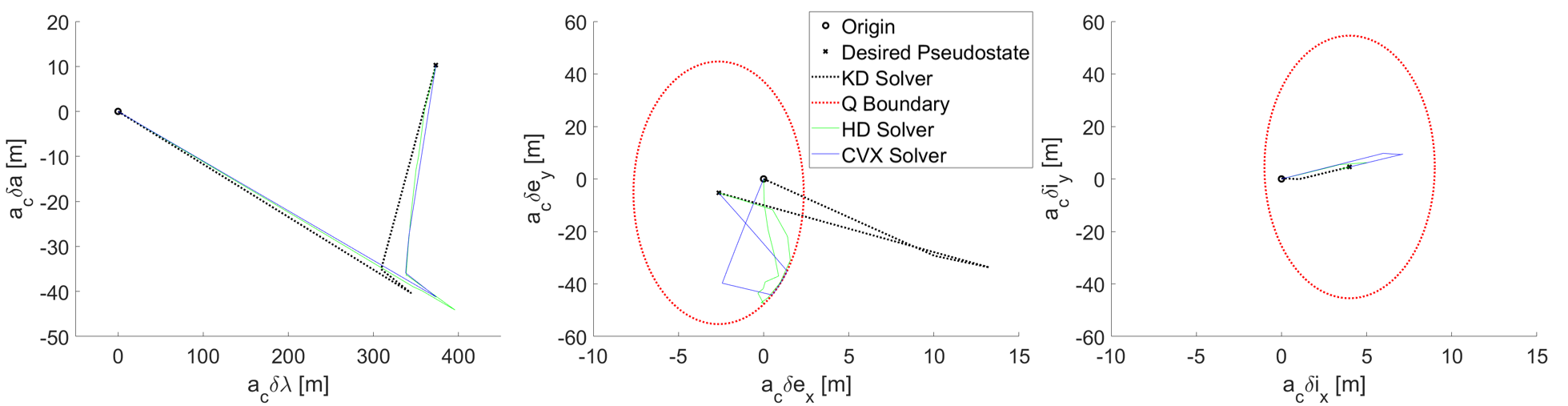}
    \caption{Pseudostate trajectory of the HD solver compared to the CVX solver for control under convex magnitude and state constraints.}
    \label{fig: convex constraints}
\end{figure*}
\subsection{Convex Constraints}

Next, the capabilities of the HD solver are demonstrated for convex state and maneuver magnitude constraints. The initial and final desired deputy states, and resulting pseudostate, that define this scenario are provided in Table \ref{tab: convex constraint roe} for an initial chief spacecraft $u = 320^\circ$ at $t_0$.
\begin{table}[b]
    \centering
    \setlength{\tabcolsep}{5.5pt}
    \caption{\bf Initial, Desired, and Pseudostate Mean ROE Deputy States for Convex Constraint Validation (m)}
    \begin{tabular}{l | c c c c c c}
         & $a_c\delta a$ & $a_c \delta \lambda$ & $a_c \delta e_x$ & $a_c \delta e_y$ & $a_c \delta i_x$ & $a_c \delta i_y$ \\
         \hline \hline
        \rule{0pt}{2ex}$\boldsymbol{x}_0$ & -10.25 & -269 & 3 & 205 & -4 & 195 \\
        $\boldsymbol{x}_f$ & 0 & 200 & 0 & 200 & 0 & 200 \\
        $\boldsymbol{\omega}$ & 10.25 & 373.27 & -2.65 & -5.26 & 4 & 4.56 \\
        \hline
    \end{tabular}
    \label{tab: convex constraint roe}
\end{table}
This scenario will consider a reconfiguration that swaps the relative locations of formation flying spacecraft, such that the deputy changes its nominal separation from the negative to the positive tangential direction through a Walking Safety Ellipse (WSE) approach. A WSE is a passively-safe formation in which the deputy uses a semi-major axis offset to change its tangential separation from the chief spacecraft through natural orbital dynamics. The in-plane (RT) and out-of-plane (RN) motion are phased $90^\circ$ apart to cause bounded, elliptical relative motion perpendicular to the direction of travel. In ROE-space, the WSE is commonly realized through a $\delta a$ offset and e-i vector separation, constraining $\delta \boldsymbol{e}$ and $\delta \boldsymbol{i}$ to be parallel or anti-parallel to set a lower bound on inter-spacecraft separation through minimum separation in the RN-plane \cite{damico_phd}. This can be translated to convex maximum deviation constraints from the nominal desired $\delta \boldsymbol{e}$ and $\delta \boldsymbol{i}$ states at the end of the reconfiguration. This scenario will define $Q$ as the set of states outside of two, 2-dimensional ellipsoids centered on the desired $\delta \boldsymbol{e}$ and $\delta \boldsymbol{i}$ pseudostates with semi-minor axes of 5m along $\delta e_x$ and $\delta i_x$ and semi-major axes of 50m along $\delta e_y$ and $\delta i_y$. This ellipsoid shape allows larger magnitudes for maneuvers that keep $\delta \boldsymbol{e}$ and $\delta \boldsymbol{i}$ parallel to preserve minimum spacecraft separation in the RN plane. Additionally, this scenario defines $\Delta c_{min}$ as 1e-4m/s, $c_{max}$ as 1e-3m/s, $\Delta t$ as 1s, and $\epsilon_{subopt}$ as 0.5 (50\%). $\Delta t$ is set at a small non-zero value to ensure that $c_{max}$ is enforced at each candidate maneuver time. Recalling that the contact function for constrained maneuvering cannot be evaluated as a single expression, Rollout instead evaluates coarse samples on the unit sphere with no angular gap between the samples greater than $60^\circ$ (Alg. \ref{alg: rollout}, Line 6). The CVX solver is configured with (\ref{eq: cvx problem}) to consider maneuver components in the R, T, and N directions at every candidate maneuver time in $T$ and constrain each waypoint $\boldsymbol{z}_i$ to remain inside the separate $\delta \boldsymbol{e}$ and $\delta \boldsymbol{i}$ 2-dimensional ellipsoids, resulting in 97523 decision variables and 62420 equality constraints. The resulting ROE trajectory of the HD solver, the CVX solver, and the KD solver are displayed in Figure \ref{fig: convex constraints}.

The HD algorithm is able to accomplish impulsive control under convex constraints, generating a pseudostate trajectory that remains inside the $\delta \boldsymbol{e}$ and $\delta \boldsymbol{i}$ bounding ellipsoids for the duration of the control window. The HD solver finds a 18-maneuver plan with a total $\Delta$v cost of 0.01073m/s, 11.74\% greater than the true optimal cost of the constrained reconfiguration found by the CVX solver as 0.009630m/s. The KD solver finds a lower optimal $\Delta$v cost of the 0.008368m/s for the unconstrained problem, corresponding to the clear constraint violation in the $\delta \boldsymbol{e}$ state constraint. The sub-optimality bound $c_{subopt}$ sums both the maneuver optimality lost to constraint satisfaction and the minimum step-sizes used to optimize the maneuver magnitudes, resulting in a bound of 0.004162m/s. The cost of the constrained maneuver plan is within $c_{subopt}$ of the cost found by the KD solver, and the ratio of $c_{subopt}$ to unconstrained optimal cost is 49.74\%, within the user-specified $\epsilon_{subopt}$ (50\%). The HD solver achieves the desired final state with a root mean squared error of 0.00007802m, while the CVX solver achieves a root mean squared error of 1.737e-11m. Both errors are $<0.01\%$ normalized to the original pseudostate. Therefore, the HD solver successfully solves impulsive control under convex constraints within the optimality tolerance set by the user. Constraint satisfaction significantly increases the required computational effort of the approach. The HD solver required a significantly shorter runtime of 229.60s compared to the 1185.98s runtime required for the CVX solver, while the KD solver reached a solution for unconstrained control in 2.00s. It is worth noting here that this type of enclosing constraint requires careful $c_{max}$ selection as coupling effects between the constraints can have detrimental effects on the HD solver. Recall that the maximum maneuver constraint $c_{max}$ dictates the removal of maneuver directions in Rollout that can enter the constraint area $Q$ (Alg. \ref{alg: rollout}, Line 6). Because the HD solver assumes worst-case constrained maneuvering, a relatively large $c_{max}$ can cause the system to be uncontrollable or the desired state to be unreachable under the remaining allowed maneuvers, and the algorithm fails to find a valid maneuver plan. Care must be taken in the selection of $c_{max}$ both in regards to the real-world actuation system and the state constraints on the dynamic system. 

The results illustrated by the above example are further verified by a Monte Carlo solver comparison for convex control. The HD solver, CVX solver, and KD solver are all used to produce maneuver plans for 200 pseudostates with a Monte Carlo control domain of $T_{MC} = \{ 0, 20, 40, ... 39000 \}$. This $T$ subset of 1951 candidate maneuver times is used to reduce the computational load of the Monte Carlo analysis. The CVX solver is configured the same as previously, resulting in 50712 decision variables and 25361 equality constraints. Each run of the Monte Carlo simulation constrains 2 ROE dimensions, chosen with uniform probability as either the $\delta a$-$\delta \lambda$, $\delta\boldsymbol{e},$ or $\delta\boldsymbol{i}$ state pairs. Each unconstrained pseudostate ROE is generated from a zero mean Gaussian distribution with 1km standard deviation, while each constrained pseudostate ROE is generated from a zero mean Gaussian distribution with 10m standard deviation. The KD solver is used before the HD solver or CVX solver to generate an unconstrained maneuver plan to inform random constraint generation that requires constrained control, such that unconstrained control cannot satisfy the constraints, while ensuring the problem remains feasibly reachable. The set $Q$ is defined as a 2-dimensional bounding ellipsoid centered at the origin of pseudostate-space, with axes along the constrained dimensions and sized with respect to the maximum distance of the unconstrained trajectory from the origin ($z_{max} = \max(||\boldsymbol{z}_i||_2)$). Both axes are sampled from Gaussian distributions with a standard deviation $0.125z_{max}$, one randomly chosen to be from mean $z_{max}$ and the other from mean $0.75z_{max}$. $c_{max}$ is randomly chosen from a uniform distribution between the minimum and maximum maneuver size of the unconstrained maneuver plan, and $\Delta c_{min}$ is set one magnitude below this value. Matching the above example convex reconfiguration, $\Delta t$ is defined as 1s, and $\epsilon_{subopt}$ is defined as 0.5 (50\%). Table \ref{tab: convex constraints} displays the minimum, mean, and maximum values associated with metrics of interest over all Monte Carlo runs. These metrics include the percentage difference between HD solver cost $c_{HD}$ and CVX solver cost $c_{CVX}$, HD solver cost relative to KD solver cost $c_{KD}$ and residual sub-optimality $c_{subopt}$, HD solver runtime, and CVX solver runtime.
\begin{table}[b]
    \centering
    \caption{\bf Important Metrics in a Monte Carlo Simulation of Convex Constrained Control Scenarios}
    \begin{tabular}{c| c c c}
        Metric & Min & Mean & Max \\
        \hline \hline
        \rule{0pt}{1ex} $(c_{HD} - c_{CVX})/c_{CVX}$ (\%) & -88.64 & -25.40 & 113.59 \\
        $c_{HD}/(c_{KD} + c_{subopt})$ (\%) & 52.27 & 64.75 & 91.25 \\
        HD Runtime (s) & 15.03 & 132.76 & 718.91 \\
        CVX Runtime (s) & 162.74 & 278.76 & 456.07 \\
        \hline
    \end{tabular}
    \label{tab: convex constraints}
\end{table}

Four metrics of the Monte Carlo analysis are important to quantify the performance of the HD solver. The percentage difference in cost between the HD solver and CVX solver quantifies the HD solver performance against the optimal solution under convex constraints. However, it is evident that there is a large difference in minimum and maximum relative cost between the solvers, and the HD solver often finds a lower-cost plan than the CVX solver. This occurs because the CVX solver often struggles to fully converge to an optimal solution, even using the reduced number of candidate maneuver times in $T_{MC}$. Therefore, the CVX solver does not always output the optimal solution, but the best case feasible solution it can achieve, and the HD solver frequently approaches the optimal solution better than the CVX solver. The ratio of HD solver cost to the sum of KD solver cost and residual optimality quantifies the expected increase in cost from the unconstrained optimal solution to the bounded sub-optimal solution under the user-defined $\epsilon_{subopt}$ and $\Delta c_{min}$. This ratio never exceeds 1 over the course of the Monte Carlo, verifying that the HD solver guarantees the sub-optimality bound $\epsilon_{subopt}$ specified by the user. Finally, the HD solver and CVX solver runtimes are used comparatively to evaluate computational efficiency. The HD solver is capable of reaching a maneuver plan significantly faster than the CVX solver for most Monte Carlo runs, up to an order of magnitude difference. Furthermore, the CVX solver runtimes are significantly shorter using $T_{MC}$ than that using $T$ in the previous illustrative convex constraint scenario, while the HD solver runtimes are comparable. This implies an exponential correlation between highly discretized state spaces and computational load in the CVX solver that the HD solver is robust against. The HD solver also experiences infrequent long-runtime outliers, 15 instances of $>$400s runtime during the Monte Carlo simulation, when the user-defined optimality and constraint inputs make the HD solver poorly conditioned for the control problem. This results in many smaller maneuvers of cost $\Delta c_{min}$, producing a large increase in overall maneuver plan cost that also corresponds to the large maximum percentage difference in cost between the HD solver and CVX solver. These instances of poor HD solver performance are a function of the Monte Carlo simulation naively randomizing inputs for control, rather than applying a more-informed knowledge of the problem in solver setup. 

\subsection{Nonconvex Constraints}
\begin{figure*}[h]
    \centering
    \includegraphics[width=1\linewidth]{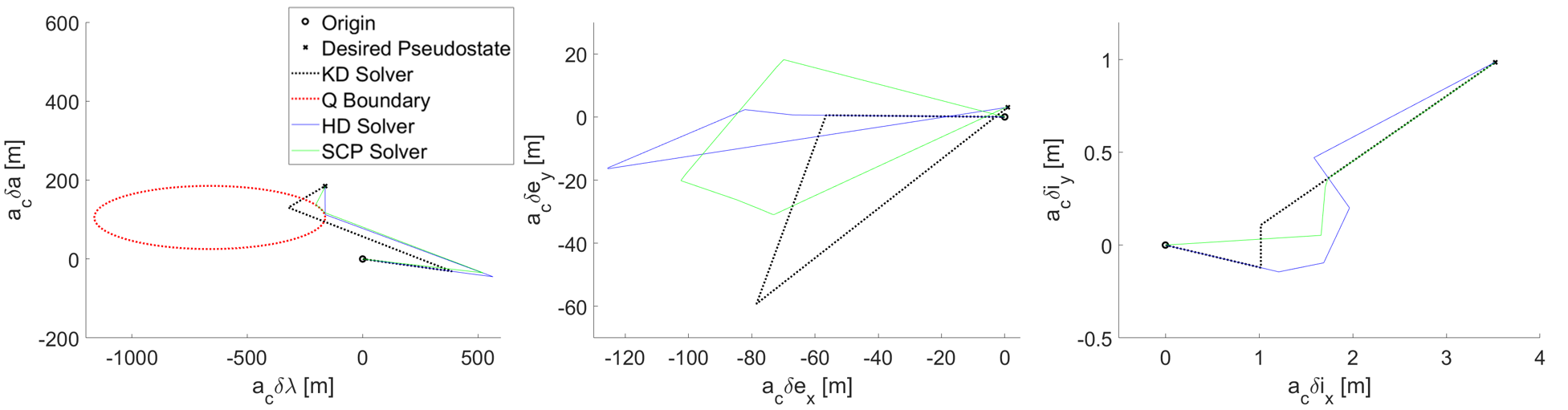}
    \caption{Pseudostate trajectory of the HD solver compared to the solution of the SCP solver for control under convex magnitude constraints and nonconvex state constraints and $\Delta t = 1$s. }
    \label{fig: nonconvex constraints dt1}
\end{figure*}
\begin{figure*}[h]
    \centering
    \includegraphics[width=1\linewidth]{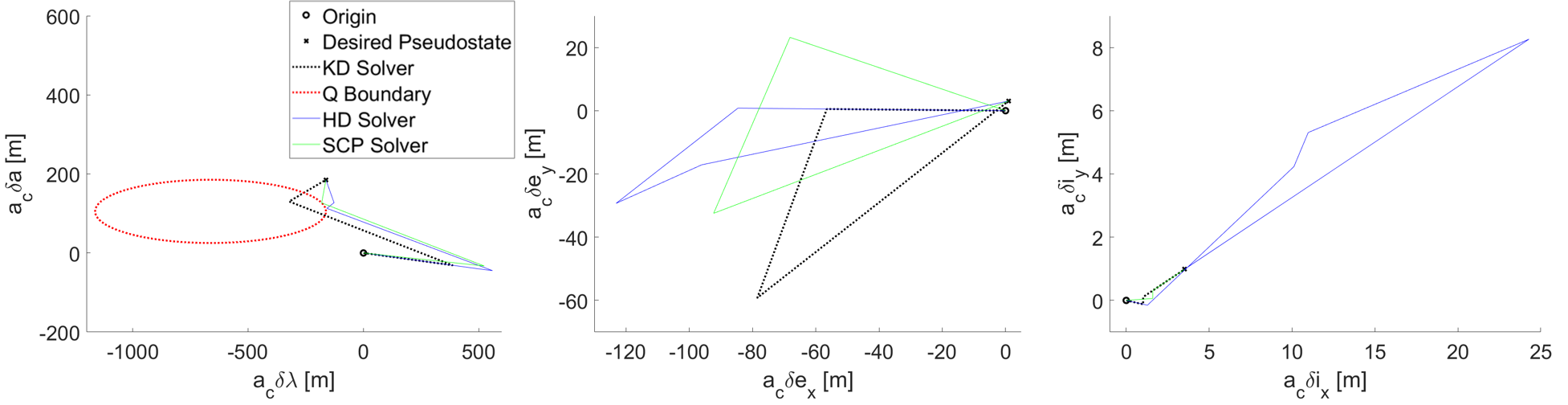}
    \caption{Pseudostate trajectory of the HD solver compared to the solution of the SCP solver for control under convex magnitude constraints and nonconvex state constraints and $\Delta t = 3000$s.}
    \label{fig: nonconvex constraints dt3000}
\end{figure*}

Finally, the capabilities of the HD solver are demonstrated for nonconvex state and maneuver timing constraints. The initial and final desired deputy states, and resulting pseudostate, that define this scenario are provided in Table \ref{tab: nonconvex constraint roe} for an initial chief spacecraft $u = 200^\circ$ at $t_0$.
\begin{table}[b]
    \centering
    \setlength{\tabcolsep}{5.5pt}
    \caption{\bf Initial, Desired, and Pseudostate Mean ROE Deputy States for Nonconvex Constraint Validation (m)}
    \begin{tabular}{l | c c c c c c}
         & $a_c\delta a$ & $a_c \delta \lambda$ & $a_c \delta e_x$ & $a_c \delta e_y$ & $a_c \delta i_x$ & $a_c \delta i_y$ \\
         \hline \hline
        \rule{0pt}{2ex}$\boldsymbol{x}_0$ & -105.25 & -321 & -1.05 & -2.73 & -3.52 & -0.51 \\
        $\boldsymbol{x}_f$ & 80 & 500 & 0 & 0 & 0 & 0 \\
        $\boldsymbol{\omega}$ & 185.25 & -162.57 & 0.93 & 3.05 & 3.52 & 0.98 \\
        \hline
    \end{tabular}
    \label{tab: nonconvex constraint roe}
\end{table}
This scenario considers a reconfiguration where the deputy approaches the chief for proximity operations in passively safe steps using a Coelliptic Orbit approach. The Coelliptic Orbit is a passively-safe formation in which the chief and deputy share the same orbital plane and use a semi-major axis offset to change tangential separation through natural dynamics with minimal periodic relative motion. In ROE-space, a Coelliptic Orbit is realized through a $\delta a$ offset and approximately zero-magnitude $\delta \boldsymbol{e}$ and $\delta \boldsymbol{i}$ vectors, and this can be represented as a KOE about the chief spacecraft in the $\delta a$ and $\delta \lambda$ states. Therefore, this scenario defines $Q$ as a two-dimensional KOE conservatively sized at the final chief state in pseudostate space with a 500m semi-major axis in the $\delta \lambda$ direction and a 80m semi-minor axis in the $\delta a$ direction, $\Delta c_{min}$ as 1e-4m/s, $c_{max}$ as 0.1m/s, and $\epsilon_{subopt}$ as 0.2 (20\%). The effect of $\Delta t$ is demonstrated comparatively by first setting $\Delta t$ as 1s, then finding a new solution for $\Delta t$ as 3000s. Rollout again evaluates coarse samples on the unit sphere with no angular gap between the samples greater than $60^\circ$ (Alg. \ref{alg: rollout}, Line 6). The approximate convex problem solved at each iteration of the SCP solver is solved by the CVX solver in a similar configuration to the previous Convex Constraint Validation, considering maneuver components in the R, T, and N directions at every 10 candidate maneuver times in $T$ ($T_{SCP} = \{0, 100, 200, ... 39000 \}$) and the constraints in (\ref{eq: scp convex approx}). This results in 391 candidate maneuver times, 12699 decision variables, and 4482 equality constraints at each SCP solver iteration. The resulting ROE trajectories of the HD solver, the SCP solver, and the KD solver are displayed in Figures \ref{fig: nonconvex constraints dt1} and \ref{fig: nonconvex constraints dt3000}.

The HD solver accomplishes impulsive control under nonconvex constraints, generating a pseudostate trajectory that remains outside of the KOE and satisfying the differing time constraints. For $\Delta t = 1$s, the HD solver finds a 5-maneuver plan with a total $\Delta$v cost of 0.03161m/s, and the sub-optimality bound $c_{subopt}$, again summing both constraint satisfaction and $\Delta c_{min}$ effects, returns as 0.004007m/s. For $\Delta t = 3000$s, the HD solver finds a 5-maneuver plan with a total $\Delta$v cost of 0.03190/s, and the sub-optimality bound $c_{subopt}$ returns as 0.005108/s. This is an expected increase in cost do to the more restrictive $\Delta t$ constraint, primarily reflected in the significant deviation of the $\delta \boldsymbol{i}$ trajectory when approaching $Q$ from the desired pseudostate with backward time-ordered maneuvers (Fig. \ref{fig: nonconvex constraints dt3000}, right). The KD solver finds a lower optimal $\Delta$v cost of the 0.02729m/s for the unconstrained problem, corresponding to the clear constraint violation of the KOE. The cost of both constrained maneuver plans is within their $c_{subopt}$ of the cost found by the KD solver, and the ratio of $c_{subopt}$ to unconstrained optimal cost is 14.68\% and 18.72\%, within the user-specified $\epsilon_{subopt}$ (20\%). The HD solver achieves the desired final state with root mean squared errors of 0.0001504m for $\Delta t = 1$s and 0.003421 for $\Delta t = 3000$s. Both errors are $<0.01\%$ normalized to the original pseudostate. Therefore, the HD solver successfully solves impulsive control under nonconvex constraints within the optimality tolerance set by the user. 

In contrast, the SCP solver finds a 78-maneuver plan with a total $\Delta$v cost of 0.02943m/s for $\Delta t = 1$s and a 5-maneuver plan with a total $\Delta$v cost of 0.02967/s for $\Delta t = 3000$s. The maneuver plans found by the HD solver are 7.41\% and 7.62\% higher than these costs, due to the restrictive assumptions used in the solver's creation. The HD solver only considers maneuvers that will not enter $Q$ from the given waypoint at the highest possible cost $c_{max}$ to prevent the solver from exploring regions of pseudostate space where the desired pseudostate becomes unreachable. This feature is clear in the $\delta a$-$\delta\lambda$ trajectory in Figure \ref{fig: nonconvex constraints dt1} for backward time-ordered maneuvers from the desired pseudostate. The SCP solver finds a more-efficient maneuver acting in the direction of $Q$ that the HD solver cannot find. However, these more-optimal costs come at the cost of numerical performance. The SCP solutions result in 1.435m root mean squared error for $\Delta t = 1$s and 0.01329m for $\Delta t = 3000$s, both significantly higher than the HD solver although still $<0.01\%$ normalized to the original pseudostate. The CVX solver used at each iteration of the SCP solver struggles to fully converge to an approximate optimal solution, resulting in this final state error. Comparing computational effort, the HD solver requires 84.65s and 76.09s to find valid constrained maneuver plans, while the SCP solver requires 471.11s and 478.24s to reach its stopping tolerance. This difference in numerical efficiency is extremely significant, recalling that the CVX solver used within the SCP solver only considers every 10 candidate maneuver times in $T$ ($T_{SCP} = \{0, 100, 200, ... 39000 \}$). Using the same control domain as the HD solver would require a similar runtime to the CVX solver in the previous Convex Constraint Validation at every iteration of the SCP solver.

The results illustrated by the above example are further verified by a Monte Carlo solver comparison for nonconvex control similar to that used for convex control. The HD solver, SCP solver, and KD solver are all used to produce maneuver plans for 200 pseudostates. The HD solver and KD solver use a Monte Carlo control domain of $T_{MC} = \{ 0, 20, 40, ... 39000 \}$. This $T$ subset of 1951 candidate maneuver times is again used to reduce the computational load of the Monte Carlo analysis. The SCP solver uses the previous control domain of $T_{SCP} = \{0, 100, 200, ... 39000 \}$ and the same internal CVX solver configuration. Each run of the Monte Carlo simulation again constrains 2 ROE dimensions, chosen with random probability as either the $\delta a$-$\delta \lambda$, $\delta\boldsymbol{e},$ or $\delta\boldsymbol{i}$ state pairs. Each pseudostate ROE is generated from a zero mean Gaussian distribution with 1000m standard deviation. The constraints are again defined with respect to the unconstrained maneuver plan from the KD solver to require constrained maneuvering while maintaining feasible reachability. The set $Q$ is defined each run as a 2-dimensional KOE on the trajectory of the unconstrained maneuver plan. The KOE center is sampled from a Gaussian with 2m standard deviation and mean at a waypoint of the unconstrained maneuver plan, chosen with uniform probability. Each KOE axis is aligned with a constrained dimension and sized from a uniform distribution between 0m and 100m. $c_{max}$ is randomly chosen from a uniform distribution between the minimum and maximum maneuver magnitude of the unconstrained maneuver plan, and $\Delta c_{min}$ is set one magnitude below this value. $\Delta t$ is sampled from a uniform distribution between 0s and 2000s. $\epsilon_{subopt}$ is defined as 0.5 (50\%) to balance constraint satisfaction and optimality. Table \ref{tab: nonconvex constraints} displays the minimum, mean, and maximum values associated with metrics of interest over all Monte Carlo runs. These metrics include the percentage difference between HD solver cost $c_{HD}$ and SCP solver cost $c_{SCP}$, HD solver cost relative to KD solver cost $c_{KD}$ and residual sub-optimality $c_{subopt}$, HD solver runtime, and SCP solver runtime.
\begin{table}[t]
    \centering
    \caption{\bf Important Metrics in a Monte Carlo Simulation of Nonconvex Constrained Control Scenarios}
    \begin{tabular}{c| c c c}
        Metric & Min & Mean & Max \\
        \hline \hline
        \rule{0pt}{1ex} $(c_{HD} - c_{SCP})/c_{SCP}$ (\%) & -0.002 & 35.99 & 182.82 \\
        $c_{HD}/(c_{KD} + c_{subopt})$ (\%) & 54.97 & 65.83 & 82.68 \\
        HD Runtime (s) & 19.15 & 95.59 & 417.63 \\
        SCP Runtime (s) & 245.26 & 301.92 & 347.63 \\
        \hline
    \end{tabular}
    \label{tab: nonconvex constraints}
\end{table}

The same four metrics are used here as the Convex Constraint Validation to assess the HD solver performance. The HD solver frequently finds comparable maneuver plan costs to the approximately optimal SCP solver, occasionally finding an equally or more optimal solution. The ratio of HD solver cost to the sum of KD solver cost and residual sub-optimality never exceeds 1 over the course of the Monte Carlo and verifies that the HD solver guarantees the sub-optimality bound $\epsilon_{subopt}$ specified by the user. Finally, the HD solver finds a constrained maneuver plan with a significantly shorter runtime than the SCP solver in almost all cases, only exceeding 200s in 17 runs of the Monte Carlo simulation. These infrequent runtime outliers again occur when the user-defined optimality and constraint inputs cause the HD solver to be poorly conditioned for the control problem. This results in many small maneuvers of cost $\Delta c_{min}$ that also cause the large outliers in cost difference between the HD solver and SCP solver. Therefore, poor performance is again a functionality of naively randomized inputs that can be mitigated by informed problem setup.

\section{Conclusion}

This work extends a reachable set theory approach for fuel-optimal impulsive control to constrained problems of interest in the field of Distributed Space Systems (DSS). A single reachable set, formed as a convex hull of all possible maneuvers over a specified control window, can be reformulated as the sequential intersection of multiple time-ordered reachable sets and broken down into a series of optimal sub-problems. This enables the derivation of necessary and sufficient conditions for not only the timing and direction of optimal maneuvering, but also the optimal magnitude of a given maneuver. Furthermore, the sequential approach enables both convex and nonconvex constraints to be implemented at the sub-problem level without violating the convexity of the reachable set. The optimality of any given maneuver can then be bounded against the unconstrained solution without the need to explicitly solve the constrained control problem. The sequential optimality conditions are derived with respect to this bound, and these new conditions are leveraged to create a constrained control solver for impulsive control that balances cost-efficiency and computational effort. The solver rolls out single maneuver constrained trajectories forward-in-time from the initial state and backward-in-time from the desired state, finding the optimal magnitude of each maneuver given the constraints of the given control problem. This process is repeated until the remaining control problem can be solved with an unconstrained control solution without violating constraints.

The overall sequential approach is demonstrated to recover the optimal cost of the unconstrained control problem within the expected optimality losses incurred by sequential maneuver magnitude optimization. The proposed solver is validated against common DSS reconfiguration scenarios, such as a Walking Safety Ellipse (WSE) or Coelliptic Orbit formation. The solver successfully finds valid maneuver plans within the user-specified optimality bound for each state, time, and maneuver magnitude constrained reconfiguration scenario with arbitrary convexity. Through Monte Carlo analysis, the proposed solver finds constrained maneuver plan costs comparable to solvers in the literature with significant improvements in algorithmic runtime and computational efficiency, especially for nonconvex problems. Overall, this work demonstrates the advantages of reachable set theory for trajectory optimization. The methodology provides insight into the nature of optimal maneuvering, solves nonconvex control problems with provably optimality, and generates control solutions with minimal computational effort.

\appendices{}

\section{Validation Dynamics Model}

This appendix describes the Linear Time Variant (LTV) motion model used in Validation (Section 6) to describe the relative motion of DSS with a quasi-nonsingular Relative Orbital Element (ROE) state, defined in Eq. \ref{eq: ROE}. The mean ROE under unperturbed Keplerian dynamics are constant values equivalent to the Integration Constants of the Hill-Clohessy-Wiltshire (HCW) and Yamanaka-Ankersen (YA) differential equations of spacecraft relative motion at small separations \cite{alfriend}. Osculating and mean ROE, similar to the Keplerian orbital elements that comprise them, vary slowly under primary orbital perturbations. Therefore, the secular and long-term evolution of ROE under J2-perturbed motion can be described by a State Transition Matrix (STM) that linearly propagates the ROE state over a given time interval, defined as
\begin{align} \label{eq: Gamma matrix}
    &\boldsymbol{\Phi}(t) \! = \! 
    \begin{bmatrix}
    1 & 0 & 0 & 0 & 0 & 0 \\
    \Phi_{21} & 1 & \Phi_{23} & \Phi_{24} & \Phi_{25} & 0 \\
    \Phi_{31} & 0 & \Phi_{33} & \Phi_{34} & \Phi_{35} & 0 \\
    \Phi_{41} & 0 & \Phi_{43} & \Phi_{44} & \Phi_{45} & 0 \\
    0 & 0 & 0 & 0 & 1 & 0 \\
    \Phi_{61} & 0 & \Phi_{63} & \Phi_{64} & \Phi_{65} & 1 \\
    \end{bmatrix}, 
\end{align}
\noindent with the simplifications
\begin{align}
    &\tau = t_f - t \text{ , } \eta_c = \sqrt{1-e_{c}^{2}} \text{ , } n_c = \sqrt{\mu_E/a_c^3}, \\
    &\kappa_c = \frac{3}{4}\frac{J_2 R^{2}_{E} \sqrt{\mu_E}}{a_{c}^{7/2} \eta_{c}^{4}} \text{ , } \Dot{\omega}_c = \kappa_c(5\cos^2(i_c)-1), \nonumber \\
    &E = 1+\eta_c \text{ , } F = 4+2\eta_c \text{ , } G = \eta^{-2} \text{ , } W = \cos^2(i_c/2), \nonumber \\
    &P = 3\cos^2(i_c)-1 \text{ , } R = \cos(i_c) \text{ , } V = \tan(i_c/2) \nonumber \\
    &e_{x1} = e_c\cos(\omega_c(t)) \text{ , } e_{y1} = e_c\sin(\omega_c(t)) \nonumber \\
    &e_{x2} = e_c\cos(\omega_c(t_f)) \text{ , } e_{y2} = e_c\sin(\omega_c(t_f)), \nonumber 
\end{align}
\noindent for the following substitutions
\begin{align}
    &\Phi_{21} = -(\frac{3}{2}n_c+\frac{7}{2}\kappa_cEP)\tau, \nonumber \\
    &\Phi_{23} = \kappa_c e_{x1}FGP\tau, \nonumber\\
    &\Phi_{24} = \kappa_c e_{y1}FGP\tau, \nonumber \\
    &\Phi_{25} = -\kappa_c F\sin(2i_c)\tau, \nonumber \\
    &\Phi_{31} = \frac{7}{2}\kappa_c e_{y2} (5\cos^2(i_c)-1)\tau, \nonumber \\
    &\Phi_{33} = \cos(\Dot{\omega}_c\tau) - 4\kappa_c e_{x1} e_{y2} G(5\cos^2(i_c)-1) \tau, \nonumber \\
    &\Phi_{34} = -\sin(\Dot{\omega}_c\tau) - 4\kappa_c e_{y1} e_{y2} G(5\cos^2(i_c)-1)\tau, \nonumber \\
    &\Phi_{35} = 5\kappa_c e_{y2} \sin(2i_c)\tau, \nonumber \\
    &\Phi_{41} = -\frac{7}{2}\kappa_c e_{x2} (5\cos^2(i_c)-1)\tau, \nonumber \\
    &\Phi_{43} = \sin(\Dot{\omega}_c\tau) + 4\kappa_c e_{x1} e_{x2} G(5\cos^2(i_c)-1) \tau, \nonumber \\
    &\Phi_{44} = \cos(\Dot{\omega}_c\tau) + 4\kappa_c e_{y1} e_{x2} G(5\cos^2(i_c)-1)\tau, \nonumber \\
    &\Phi_{45} = -5\kappa_c e_{x2} \sin(2i_c)\tau, \nonumber \\
    &\Phi_{61} = \frac{7}{2}\kappa_c\sin(2i_c)\tau, \nonumber \\
    &\Phi_{63} = -4\kappa_c e_{x1} G\sin(2i_c)\tau, \nonumber \\
    &\Phi_{64} = -4\kappa_c e_{y1} G\sin(2i_c)\tau, \nonumber \\
    &\Phi_{65} = 2\kappa_c \sin^2(i_c)\tau, \nonumber
\end{align}
\noindent where $\boldsymbol{\Phi}(t) = \boldsymbol{\Phi} (\boldsymbol{\alpha}_c(t),t_f-t)$, $\boldsymbol{\alpha}_c(t)$ are the mean Keplerian orbital elements ($a$, $e$, $i$, $\Omega$, $\omega$, and $M$) at time $t$, subscript $c$ denotes the chief spacecraft, $t_f$ is the end of the control window, $n$ is the mean motion, $R_E$, $J_2$, and $\mu_E$ are the radius, oblateness, and standard gravitational parameter of Earth, and $\Phi_{ij}$ is the element in the $i$th row of the $j$th column \cite{koenig_STM}. RTN delta-v control actions can be mapped to changes in osculating ROE, and approximately to changes in mean ROE, through the Gauss Variational Equations (GVEs). This mapping is described by control input matrix $\boldsymbol{B}$ for orbits of arbitrary eccentricity as
\begin{equation} \label{eq: B matrix}
    \hspace{-2mm}
    \boldsymbol{B}(t) = 
    \frac{1}{a_c n_c}
    \begin{bmatrix}
    B_{11} \! \! & B_{12} \! \! & 0 \\
    B_{21} \! \! & 0 \! \! & 0 \\
    B_{31} \! \! & B_{32} \! \! & B_{33} \\
    B_{41} \! \! & B_{42} \! \! & B_{43} \\
    0 \! \! & 0 \! \! & B_{53} \\
    0 \! \! & 0 \! \! & B_{63} \\
    \end{bmatrix},
\end{equation}
\noindent with the substitutions
\begin{align}
    &B_{11} = \frac{2}{\eta_c} e_c \sin(\nu_c) \text{ , } B_{12} = \frac{2}{\eta_c} (1 + e_c \cos(\nu_c)), \nonumber \\
    &B_{21} = -\frac{2\eta_c^2}{1 + e_c \cos(\nu_c)} \text{ , } B_{31} = \eta_c \sin(\theta_c) \nonumber \\
    &B_{32} = \eta_c \frac{(2 + e_c \cos(\nu_c))\cos(\theta_c) + e_c \cos(\omega_c)}{1 + e_c\cos(\nu_c)}, \nonumber \\
    &B_{33} = \frac{\eta_c e_c \sin(\omega_c) \sin(\theta_c)}{\tan(i_c) (1 + e_c \cos(\nu_c))} \text{ , } B_{41} = -\eta_c \cos(\theta_c), \nonumber \\
    &B_{42} = \eta_c \frac{(2 + e_c \cos(\nu_c))\sin(\theta_c) + e_c \sin(\omega_c)}{1 + e_c\cos(\nu_c)}, \nonumber \\
    &B_{43} = -\frac{\eta_c e_c \cos(\omega_c) \sin(\theta_c)}{\tan(i_c) (1 + e_c \cos(\nu_c))}, \nonumber \\
    &B_{53} = \frac{\eta_c \cos(\theta_c)}{1 + e_c \cos(\nu_c)} \text{ , } B_{63} = \frac{\eta_c \sin(\theta_c)}{1 + e_c \cos(\nu_c)}, \nonumber
\end{align}
\noindent where $\boldsymbol{B}(t) = \boldsymbol{B}(\boldsymbol{\alpha}_c(t))$, $\theta = \omega + \nu$  and $\nu$ is the true anomaly, and $B_{ij}$ is the element in the $i$th row of the $j$th column \cite{koenig_RST}. Previous authors have demonstrated that this LTV dynamic propagation model tracks the numerical integration of the nonlinear, J2-perturbed GVEs at sub-meter-level accuracy for inter-spacecraft separations under 700m and $e_c \leq 0.5$ and at meter-level accuracy for inter-spacecraft separations up to 1km and $e_c \leq 0.7$ \cite{sullivan}.

\acknowledgments
This work is supported by Blue Origin (SPO \#299266) as an Associate Member and Co-Founder of the Stanford’s Center of AEroSpace Autonomy Research (CAESAR). This article solely reflects the opinions and conclusions of its authors and not any Blue Origin entity.

\bibliographystyle{IEEEtran}
\bibliography{references.bib}

\begin{biographywithpic}
{Matthew Hunter}{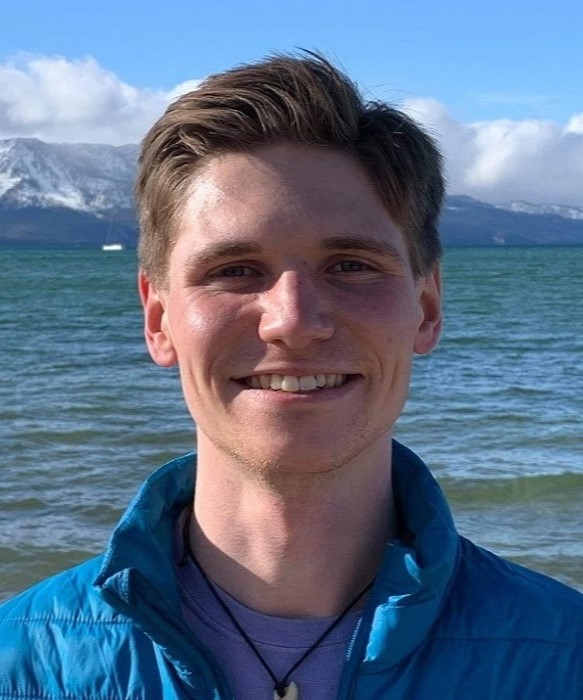}
is a Ph.D. student in the Space Rendezvous Laboratory. He “got out” of Georgia Tech with a Bachelor of Science in Mechanical Engineering, where he received the George W. Woodruff School of Mechanical Engineering Outstanding Scholar Award. He has interned at a variety of companies in the space and defense industry, including NASA Goddard, Raytheon Technologies, and Bell. Matthew is currently developing guidance and control algorithms for the Space Weather Atmospheric Multiscale EXperiment (SWARM-EX), and his research focuses on computationally efficient, optimal hybrid control techniques that incorporate perturbations to Keplerian orbit into propulsive maneuver schemes to decrease their delta-v cost. This work will enable fully integrated roto-translational spacecraft control to increase mission lifetime and enhance maneuver accuracy. 
\end{biographywithpic}

\begin{biographywithpic}
{Simone D’Amico}{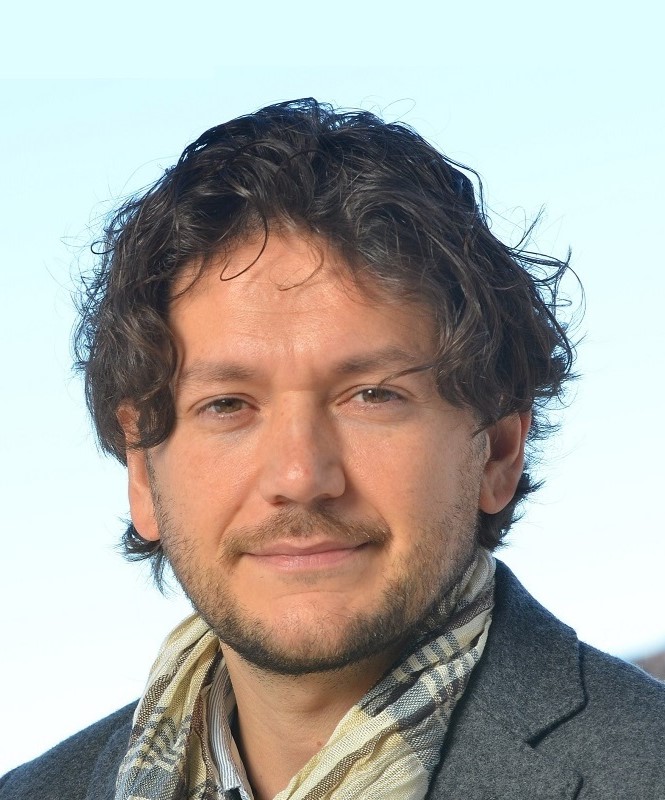}
Simone D’Amico is Associate Professor of Aeronautics and Astronautics (AA), W.M. Keck Faculty Scholar in the School of Engineering, and Professor of Geophysics (by Courtesy). He is the Founding Director of the Stanford Space Rendezvous Laboratory, Co-Director of the Center for AEroSpace Autonomy Research (CAESAR), and Director of the Undergraduate Program in Aerospace Engineering at Stanford. He has 20+ years of experience in research and development of autonomous spacecraft and distributed space systems. He developed the distributed Guidance, Navigation, and Control (GNC) system of several formation-flying and rendezvous missions and is currently the institutional PI of four autonomous satellite swarms funded by NASA (STARLING, STARI) and by NSF (VISORS, SWARM-EX). Besides academia, Dr. D’Amico is in the Advisory Board of four space start-ups focusing on distributed space systems for future applications in SAR remote sensing, orbital lifetime prolongation, and space-based solar power. He was the recipient of several awards, most recently the 2024 NASA Ames Honor Award for the Starling mission, Best Paper Awards at IAF (2022), IEEE (2021), AIAA (2021), AAS (2019) conferences, and the M. Barry Carlton Award by IEEE (2020). He received the B.S. and M.S. degrees from Politecnico di Milano (2003) and the Ph.D. degree from Delft University of Technology (2010).
\end{biographywithpic}

\end{document}